\pdfoutput=1
\documentclass[11pt]{article}
\usepackage[T1]{fontenc}
\usepackage{lmodern}

\usepackage[margin=1in]{geometry}
\usepackage{amsmath,amssymb,amsthm,mathtools}
\usepackage{microtype}
\usepackage{graphicx}
\usepackage{tikz}
\usepackage[hidelinks]{hyperref}
\usepackage{color}
\hypersetup{pdftitle={The One-Skeleton of Multilinear Polytopes and Their Standard Relaxations},pdfauthor={Alberto Del Pia, Zihong Xu}}
\newtheorem{theorem}{Theorem}
\newtheorem{lemma}{Lemma}
\newtheorem{proposition}{Proposition}
\newtheorem{corollary}{Corollary}
\theoremstyle{definition}
\newtheorem{definition}{Definition}
\theoremstyle{remark}
\newtheorem{remark}{Remark}
\theoremstyle{plain}
\DeclareMathOperator{\MP}{MP}
\DeclareMathOperator{\SR}{SR}

\DeclareMathOperator{\co}{co}
\DeclareMathOperator{\dist}{dist}
\DeclareMathOperator{\diam}{diam}
\DeclareMathOperator{\conv}{conv}
\DeclareMathOperator{\piv}{piv}

\title{The One-Skeleton of Multilinear Polytopes\\ and Their Standard Relaxations}
\author{Alberto Del Pia\thanks{Department of Industrial and Systems Engineering \& Wisconsin Institute for Discovery, University of Wisconsin--Madison. Email: \texttt{delpia@wisc.edu}.}
\and Zihong Xu\thanks{Department of Mathematics, University of Wisconsin--Madison. Email: \texttt{zxu455@wisc.edu}.}}
\date{November 11, 2026}

\begin{document}
\maketitle
\begin{abstract}
Edge-following methods, such as the simplex method, move along edges of a polytope. Their behavior depends on which vertices are adjacent and on the diameter. We study both for the two polytopes of binary polynomial optimization: the multilinear polytope and its standard linear relaxation. For quadratic polynomials, Padberg showed that two integral vertices are adjacent in both polytopes if and only if the variables where they differ induce a connected subgraph. For higher degree, we replace this subgraph by a contracted hypergraph on the same variables, built from the monomials whose constraints tie these variables together. We show that its connectivity decides adjacency in the relaxation but, unlike in the quadratic case, not always in the multilinear polytope. Adjacency is the same in both polytopes exactly when the relaxation has the Trubin property. We investigate when this property holds and characterize it for cubic polynomials. The contracted hypergraph also yields short paths. The diameter of the multilinear polytope is at most the number of variables, and that of the relaxation is linear in the numbers of variables and monomials. Similar linear bounds hold for the number of pivots between bases of the relaxation.
\end{abstract}

\noindent\textbf{Keywords:} Multilinear polytope, Standard relaxation, Adjacency, Trubin property, Diameter, Pivots

\section{Introduction}
Binary polynomial optimization asks to maximize a polynomial over binary variables. Its standard linearization leads to two polytopes: the multilinear polytope, whose vertices are the feasible solutions, and its standard linear relaxation. Edge-following methods, such as the simplex method and local search, move between adjacent vertices of a polytope. For such methods, adjacency determines which moves are available, and the diameter is a lower bound on the number of moves needed in the worst case. For quadratic polynomials, Padberg \cite{Pad89} showed that a single connectivity condition decides adjacency of integral vertices in both polytopes. In this paper we study adjacency and the diameter of both polytopes for polynomials of higher degree.

\paragraph{Problem setup.} We now fix notation for \emph{binary polynomial optimization} (BPO). Since $x^k=x$ for every $x\in\{0,1\}$ and every positive integer $k$, we may assume that the objective polynomial is \emph{multilinear}, that is, each variable appears in it with degree at most one. Following \cite{dPKha17MOR}, we encode the polynomial by a hypergraph. A \emph{hypergraph} $H$ is a pair $(V,E)$, where $V$ is a finite set of \emph{nodes} and $E$ is a set of subsets of $V$ of cardinality at least two, called \emph{edges} of $H$. Each binary variable $x_v$ is represented by a node $v\in V$, and each monomial $\prod_{v\in e}x_v$ of degree at least two by an edge $e\in E$. We linearize the problem by writing $z_v$ for $x_v$ and by introducing an auxiliary variable $z_e=\prod_{v\in e}z_v$ for every edge $e\in E$. The convex hull of the feasible points of the resulting problem is the \emph{multilinear polytope}
\[
    \MP(H)=\conv\Bigl\{z\in\{0,1\}^{V\cup E}: z_e=\prod_{v\in e}z_v \ \ \forall e\in E\Bigr\}.
\]
Every BPO problem whose monomials are the edges of $H$ is equivalent to the maximization of a linear function over $\MP(H)$.
For every $S\subseteq V$, the \emph{multilinear point} $z^S\in \{0,1\}^{V\cup E}$ is defined by $z^S_v=1$ if $v\in S$ and $z^S_v=0$ otherwise, for $v\in V$, and by $z^S_e=\prod_{v\in e}z^S_v$ for $e\in E$. Thus $z^S_e = 1$ if $e\subseteq S$, and $z^S_e = 0$ otherwise. The vertices of $\MP(H)$ are precisely the multilinear points $z^S$ with $S\subseteq V$, and we also call them \emph{multilinear vertices}. The \emph{standard relaxation} $\SR(H)$ of $\MP(H)$ is obtained by replacing each equation $z_e=\prod_{v\in e}z_v$ with the linear inequalities that describe the convex hull of its binary solutions \cite{Cra93}, and by replacing the integrality of the node variables with the bounds $0\le z_v\le 1$:
\[
    \SR(H) = \left\{z\in \mathbb{R}^{V\cup E}\middle|
    \begin{array}{ll}\ 0\le z_v\leq 1 & \forall v\in V,\\
    z_e\geq 0 &\forall e\in E,\\
    z_e\geq \textstyle\sum_{v\in e} z_v - |e| + 1 &\forall e\in E,\\
    z_e\leq z_v & \forall e\in E, \forall v\in e
    \end{array}\right\}.
\]
We have $\MP(H)\subseteq \SR(H)$, and the integral points of $\SR(H)$ are exactly the multilinear vertices of $\MP(H)$. Since $\SR(H)\subseteq[0,1]^{V\cup E}$, these points are also vertices of $\SR(H)$.

\paragraph{Definitions.} We write $A\triangle B$ for the symmetric difference of two sets $A$ and $B$. The \emph{rank} of a hypergraph $H$ is the largest cardinality of an edge of $H$, and a hypergraph whose edges all have cardinality two is a graph. A \emph{path} in $H$ is a sequence $(v_1,e_1,v_2,\dots,e_k,v_{k+1})$ of distinct nodes and distinct edges of $H$ such that $v_i,v_{i+1}\in e_i$ for every $i$. A \emph{connected component} of $H$ is a maximal nonempty set of nodes any two of which are joined by a path, and $H$ is \emph{connected} if it has exactly one connected component. In particular, a node contained in no edge forms a connected component by itself. Two vertices of a polytope $P$ are \emph{adjacent} if they are the endpoints of a one-dimensional face of $P$, and the \emph{one-skeleton} of $P$ is the graph on the vertices of $P$ whose edges join adjacent vertices. The \emph{diameter} $\diam(P)$ is the maximum, over all pairs of vertices of $P$, of the minimum number of edges of a path between them in the one-skeleton.

\paragraph{The graph case.} When $H$ is a graph $G$, the polytope $\MP(G)$ is the \emph{Boolean quadric polytope} of $G$. Its facial structure and adjacency properties, together with the fractional vertices of its standard relaxation, were studied in depth by Padberg \cite{Pad89}. One of his results shows that a single connectivity condition decides adjacency both in the Boolean quadric polytope and in its standard relaxation.
\begin{theorem}[{\cite[Theorem~8 and Remark~2]{Pad89}}]\label{thm:connectivity-in-graph}
    Let $G$ be a graph and let $S,T\subseteq V$ be distinct. Then the following statements are equivalent:
    \begin{itemize}
        \item[(i)] the subgraph of $G$ induced by $S \triangle T$ is connected,
        \item[(ii)] $z^S$ and $z^T$ are adjacent in $\MP(G)$,
        \item[(iii)] $z^S$ and $z^T$ are adjacent in $\SR(G)$.
    \end{itemize}
\end{theorem}
\noindent Given two polytopes $P\subseteq Q$ such that every vertex of $P$ is also a vertex of $Q$, we say that $Q$ has the \emph{Trubin property} with respect to $P$ if the one-skeleton of $P$ is contained in the one-skeleton of $Q$. The name refers to Trubin \cite{Tru69}, who established this property for the set partitioning polytope and its linear relaxation. Balas and Padberg used it to show that set partitioning problems can in principle be solved by a variant of the simplex method that visits only integral vertices (see the survey \cite{BalPad76}). The equivalence of (ii) and (iii) says that $\SR(G)$ has the Trubin property with respect to $\MP(G)$.

\begin{figure}[t]
\centering
\begin{tikzpicture}[scale=1.0,
  snode/.style={circle,fill=black,inner sep=1.6pt},
  tnode/.style={circle,draw=black,fill=white,inner sep=1.6pt}]
  \begin{scope}[shift={(0,0.35)}]
    \draw[rounded corners=10pt,fill=gray!18,draw=gray!70] (-0.55,-0.45)--(1.75,-0.45)--(0.60,1.70)--cycle;
    \node[snode] (a1) at (0,0) {}; \node[snode] (a2) at (1.2,0) {}; \node[tnode] (a3) at (0.6,0.95) {};
    \draw[thick] (a1)--(a2);
    \node[below=1pt] at (a1) {\small $1$}; \node[below=1pt] at (a2) {\small $2$}; \node[right=2pt] at (a3) {\small $3$};
    \node at (0.6,-0.85) {$H$};
  \end{scope}
  \begin{scope}[shift={(3.1,0.35)}]
    \node[snode] (b1) at (0,0) {}; \node[snode] (b2) at (1.2,0) {}; \node[tnode] (b3) at (0.6,0.95) {};
    \draw[thick] (b1)--(b2);
    \node[below=1pt] at (b1) {\small $1$}; \node[below=1pt] at (b2) {\small $2$}; \node[right=2pt] at (b3) {\small $3$};
    \node at (0.6,-0.85) {$\Theta(H,S,T)$};
  \end{scope}
  \node at (2.15,-1.1) {(a)};
  \begin{scope}[shift={(7.2,-0.45)}]
    \draw[->] (-.16,0)--(2.30,0) node[right] {$t$};
    \draw[->] (0,-.14)--(0,2.30) node[above] {$s$};
    \coordinate (O) at (0,0);\coordinate (S) at (2,0);
    \coordinate (F) at (1,2);\coordinate (T) at (0,2);
    \draw[thick] (O)--(S)--(F)--(T)--cycle;
    \draw[dashed] (S)--(T);
    \fill (O) circle (1pt) node[below left] {$z^\emptyset$};
    \fill (S) circle (1pt) node[below] {$z^S$};
    \fill (T) circle (1pt) node[left] {$z^T$};
    \fill (F) circle (1pt) node[above right] {$\bar z$};
    \node at (1.0,-0.65) {(b)};
  \end{scope}
\end{tikzpicture}
\caption{(a) The hypergraph $H$ with the edges $\{1,2\}$ and $\{1,2,3\}$ (shaded), for $S=\{1,2\}$ (black nodes) and $T=\{3\}$ (white node), and the contracted hypergraph $\Theta(H,S,T)$, in which node $3$ is isolated. (b) The minimal face of $\SR(H)$ containing $z^S$ and $z^T$, in the coordinates $t=z_1=z_2=z_{12}$ and $s=z_3$. The vertex $\bar z=(\frac12,\frac12,1,\frac12,0)$ is fractional, and the dashed segment $[z^S,z^T]$ is not a one-dimensional face.}\label{fig:counterexample}
\end{figure}
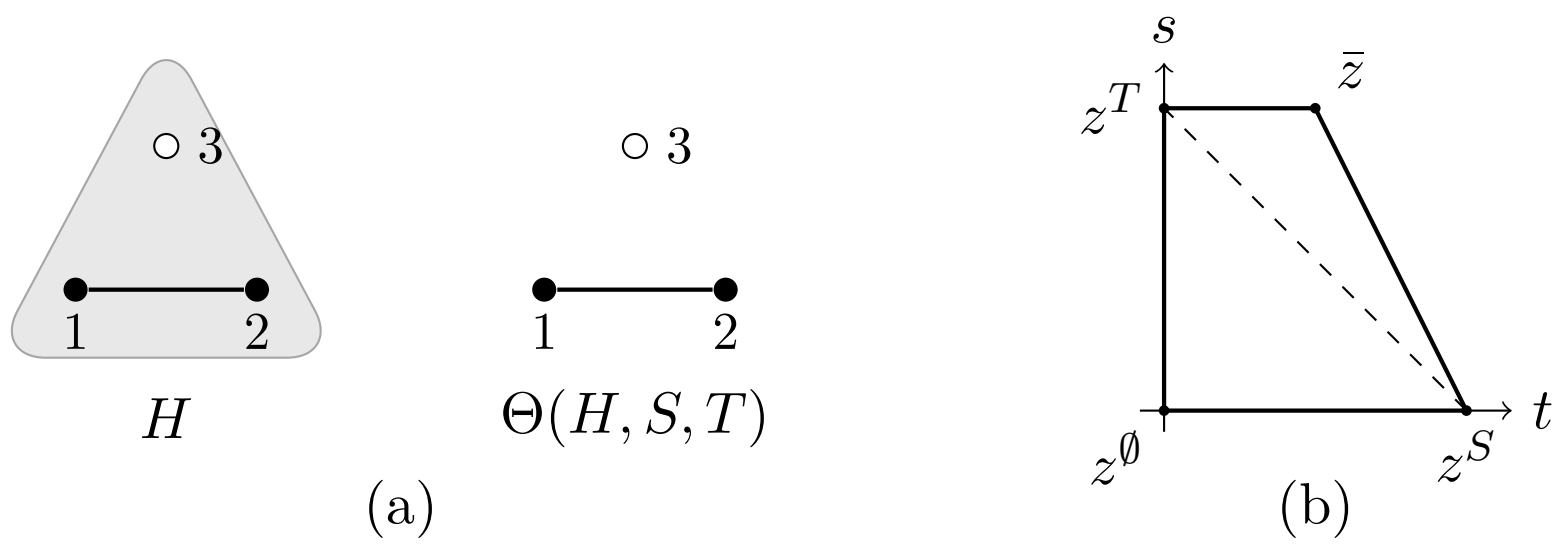

\paragraph{The contracted hypergraph.} For hypergraphs of rank at least three, to the best of our knowledge, adjacency of multilinear vertices had not been characterized in either polytope. Theorem~\ref{thm:connectivity-in-graph} suggests the question that guides our study of adjacency:
\begin{center}
    \textit{What takes over the role of ``the subgraph induced by $S\triangle T$'' when $H$ is a hypergraph?}
\end{center}
Our answer is a contracted hypergraph $\Theta(H,S,T)$ on the node set $S\triangle T$, defined as follows.
\begin{definition}\label{def:contracted-hypergraph}
    Let $H$ be a hypergraph and let $S,T\subseteq V$. For an edge $f\subseteq S\cup T$, we call $f\cap(S\triangle T)$ the \emph{trace} of $f$. We denote by $F_\Theta$ the set of all edges $f\in E$ with $f\subseteq S\cup T$ that satisfy one of the following three conditions:
    \begin{itemize}
        \item[(a)] $|f \cap (S\setminus T)| \ge 2$ and $|f \cap (T\setminus S)| = 0$,
        \item[(b)] $|f \cap (S\setminus T)| = 0$ and $|f \cap (T\setminus S)| \ge 2$,
        \item[(c)] $|f \cap (S\setminus T)| = 1$ and $|f \cap (T\setminus S)| = 1$.
    \end{itemize}
    The edges in $F_\Theta$ are the \emph{generating edges}. The \emph{contracted hypergraph} $\Theta(H,S,T)$ is the hypergraph with node set $S\triangle T$ whose edges are the traces of the generating edges, that is, with edge set
    \[
        E_\Theta := \{f\cap (S\triangle T): f\in F_\Theta\}.
    \]
\end{definition}
\noindent Intuitively, the generating edges are the edges whose constraints of $\SR(H)$, tight at the midpoint of $z^S$ and $z^T$, tie two nodes of $S\triangle T$ together, so that these nodes cannot move independently within the minimal face containing $z^S$ and $z^T$. For example, let $H$ have the edges $\{1,2\}$ and $\{1,2,3\}$, and let $S=\{1,2\}$ and $T=\{3\}$. The edge $\{1,2\}$ satisfies (a). The edge $\{1,2,3\}$ meets $S\setminus T$ in two nodes and $T\setminus S$ in one, so it is not a generating edge. Hence $\Theta(H,S,T)$ has the single edge $\{1,2\}$, and node $3$ is isolated, although an edge of $H$ contains all three nodes (Fig.~\ref{fig:counterexample}(a)).

Our first result shows that the contracted hypergraph decides adjacency in $\SR(H)$.
\begin{theorem}\label{thm:adjacency-sr}
    Let $H$ be a hypergraph and let $S,T\subseteq V$ be distinct. Then the multilinear vertices $z^S$ and $z^T$ are adjacent in $\SR(H)$ if and only if $\Theta(H,S,T)$ is connected.
\end{theorem}
\noindent When $H$ is a graph, an edge of $H$ is a generating edge precisely when both of its endpoints lie in $S\triangle T$. Hence $\Theta(H,S,T)$ is the subgraph of $H$ induced by $S\triangle T$, and Theorem~\ref{thm:adjacency-sr} recovers the equivalence of (i) and (iii) in Theorem~\ref{thm:connectivity-in-graph}. Table~\ref{tab:summary} summarizes our results and compares them with the graph case.

\begin{table}[t]
\centering
\caption{The graph case and the hypergraph case. Here $G[S\triangle T]$ denotes the subgraph of $G$ induced by $S\triangle T$. The entries for graphs follow from Padberg's results (Theorem~\ref{thm:connectivity-in-graph} and \cite[Theorem~7 and p.~157]{Pad89}), except the diameter of $\SR(G)$ (Theorem~\ref{thm:diameter-sr-graph}) and the pivot bounds, which are new.}\label{tab:summary}
\small\renewcommand{\arraystretch}{1.15}\setlength{\tabcolsep}{10pt}
\resizebox{\textwidth}{!}{%
\begin{tabular}{@{}lll@{}}
\hline
 & Graphs $G$ & Hypergraphs $H$\\ \hline
Adjacency in $\SR$ & $G[S\triangle T]$ connected & $\Theta(H,S,T)$ connected (Thm.~\ref{thm:adjacency-sr})\\
Adjacency in $\MP$ & $G[S\triangle T]$ connected & coNP-complete (Prop.~\ref{prop:mp-hardness})\\
Trubin property & always & iff no bad pair, if rank $\le3$ (Thm.~\ref{thm:rank-three})\\
Vertices of $\SR$ & half-integral & large denominators (Prop.~\ref{prop:denominators})\\
Diameter of $\MP$ & $\le|V|$ & $\le|V|$ (Thm.~\ref{thm:diameter-mp})\\
Diameter of $\SR$ & $\le|V|$ & $\le|V|+2|E|+2\min\{|V|,|E|\}$ (Thm.~\ref{thm:diameter-sr})\\
Pivots, chosen bases & $\le|V|$ & $\le|V|$ (Thm.~\ref{thm:pivot-flip})\\
Pivots, any bases & $\le3|V|+4|E|+2\min\{|V|,|E|\}$ & $\le3|V|+4|E|+2\min\{|V|,|E|\}$ (Thm.~\ref{thm:pivot-sr})\\ \hline
\end{tabular}}
\end{table}

\paragraph{Adjacency.} Theorem~\ref{thm:adjacency-sr} shows that adjacency in $\SR(H)$ remains a connectivity question. In contrast, building on a result of Matsui \cite{Mat95}, we show that deciding whether two multilinear vertices are adjacent in $\MP(H)$ is coNP-complete, even if every edge of $H$ has cardinality three (Proposition~\ref{prop:mp-hardness}). Since $\MP(H)\subseteq\SR(H)$, adjacency in $\SR(H)$ between two multilinear vertices implies adjacency in $\MP(H)$ (Lemma~\ref{lem:adjacency-subpolytope} in the appendix). Therefore the two one-skeletons agree on the multilinear vertices exactly when $\SR(H)$ has the Trubin property with respect to $\MP(H)$. A vertex of a polytope maximizes a linear function if and only if no adjacent vertex has a larger value. Hence, under the Trubin property, a multilinear vertex $z^S$ maximizes a linear function over $\MP(H)$ if and only if no $z^T$ with $\Theta(H,S,T)$ connected has a larger value.

\paragraph{The Trubin property.} By Theorem~\ref{thm:adjacency-sr}, $\SR(H)$ has the Trubin property with respect to $\MP(H)$ if and only if, for all distinct $S,T\subseteq V$,
\[
    z^S \text{ and } z^T \text{ are adjacent in } \MP(H) \quad\Longrightarrow\quad \Theta(H,S,T) \text{ is connected}.
\]
This implication fails already for a hypergraph with three nodes and two edges, the smallest possible size (Proposition~\ref{prop:sparse-counterexample} and Remark~\ref{rem:minimality} in the appendix). Our main results on the Trubin property locate the obstruction. We call a pair $S,T$ \emph{bad} if the trace of some edge $e\subseteq S\cup T$ meets exactly two connected components of $\Theta(H,S,T)$. If $H$ has a bad pair, then $\SR(H)$ does not have the Trubin property (Theorem~\ref{thm:bad-pair}). For cubic polynomials, that is, for hypergraphs of rank at most three, this obstruction is the only one (Theorem~\ref{thm:rank-three}). Hence, for these hypergraphs, a failure of the Trubin property has a certificate that can be checked in polynomial time. For hypergraphs of arbitrary rank, we give two sufficient conditions. The first one concerns a single pair $S,T$. The implication above holds for this pair if the pair is not bad and the incidence graph between the components and the traces that meet at least two of them is a forest (Theorem~\ref{thm:acyclic-traces}). The second one concerns $H$ alone: for every edge $e$, the pairs of nodes of $e$ that are not themselves edges of $H$ form a matching (Theorem~\ref{thm:matching-trubin}). This condition holds for all graphs and allows edges of arbitrary cardinality. Geometrically, every failure of the Trubin property produces a fractional vertex of $\SR(H)$ in the minimal face containing $z^S$ and $z^T$. Under a condition $(\star)$ on $S,T$ that holds for all graphs, this face is instead a cube with multilinear vertices (Theorem~\ref{thm:cube}).

\paragraph{Diameter.} The contracted hypergraph also yields short paths. \emph{Switching} a set $C\subseteq V$ at $z^S$ means moving from $z^S$ to $z^{S\triangle C}$. Switching the connected components of $\Theta(H,S,T)$ one at a time gives a path from $z^S$ to $z^T$ in the one-skeletons of both polytopes whose length is the number of connected components of $\Theta(H,S,T)$ (Theorem~\ref{thm:diameter-mp}). In particular, $\diam(\MP(H))\le|V|$, which for graphs also follows from Theorem~\ref{thm:connectivity-in-graph}. Our main diameter results concern the relaxation. Padberg proved that the vertices of $\SR(G)$ are $\{0,\frac12,1\}$-valued \cite[Theorem~7 and p.~157]{Pad89}. Using this, we show that $\diam(\SR(G))\le|V|$ (Theorem~\ref{thm:diameter-sr-graph}). Already for cubic polynomials, the vertices of $\SR(H)$ can have exponentially large denominators, so the known diameter bounds in terms of the denominators \cite{KleOnn92,dPMic16DCG} can be exponential in $|V|$ (Proposition~\ref{prop:denominators}). We work instead with bases, that is, sets of $|V|+|E|$ linearly independent constraints of $\SR(H)$ tight at a vertex. With them we prove that $\diam(\SR(H))\le|V|+2|E|+2\min\{|V|,|E|\}$, which is linear in the dimension of $\SR(H)$ (Theorem~\ref{thm:diameter-sr}).

\paragraph{Pivots.} Bases also lead to pivots, the steps of the simplex method, which are defined in Section~\ref{sec:diameter}. Since a degenerate pivot changes the basis without moving to another vertex, the diameter bounds do not directly bound the number of pivots. Between two multilinear vertices, suitably chosen bases are at most $|V|$ pivots apart, along bases at multilinear vertices only (Theorem~\ref{thm:pivot-flip}). If the initial basis is prescribed, $|E|+1=2^{|V|-1}$ pivots may be needed even between adjacent vertices (Proposition~\ref{prop:lowerbound}). Hence, for prescribed bases, no bound polynomial in $|V|$ alone exists. Between prescribed bases at arbitrary vertices of $\SR(H)$, at most $3|V|+4|E|+2\min\{|V|,|E|\}$ pivots suffice (Theorem~\ref{thm:pivot-sr}).

\paragraph{Related work.} The facial structure of the multilinear polytope has been studied extensively, starting with \cite{dPKha17MOR}. The standard relaxation coincides with $\MP(H)$ if and only if $H$ is Berge-acyclic \cite{dPKha18SIOPT,BucCraRod19}. For the larger class of kite-free $\beta$-acyclic hypergraphs, a description of $\MP(H)$ is obtained in \cite{dPKha21MOR} by means of further inequalities. Extended formulations of $\MP(H)$ are known for $\beta$-acyclic hypergraphs \cite{dPKha24MPA}, hypergraphs of bounded nest-set gap \cite{dPKha24xx}, $\alpha$-acyclic hypergraphs of bounded rank \cite{dPKha26xMPA}, and hypergraphs of bounded incidence treewidth \cite{CapdPDiG26xMPA}. For the cut polytope, which is linearly isomorphic to the Boolean quadric polytope \cite{DeS90}, Barahona and Mahjoub characterized adjacency \cite{BarMah86}. Deciding non-adjacency of two vertices is NP-complete for several $0/1$ polytopes, starting with the traveling salesman polytope \cite{Pap78}. Proposition~\ref{prop:mp-hardness} relies on a result of this type by Matsui \cite{Mat95}. For the fractional matching polytope, another relaxation with half-integral vertices, Sanit\`a characterized the diameter exactly \cite{San18}. The fractional stable set polytope is affinely isomorphic to the face of $\SR(G)$ on which all edge coordinates vanish. Using a graph characterization of its bases and pivots, Michini and Sassano proved that its diameter is at most $|V|$ \cite{MicSas14}. For a different extension of the Boolean quadric polytope, associated with satisfiability problems, Nikolaev \cite{Nik18} studied the one-skeleton and the denominators of the fractional vertices of its LP relaxation.

\paragraph{Organization.} Section~\ref{sec:adjacency} outlines the proof of Theorem~\ref{thm:adjacency-sr} and shows that adjacency in $\MP(H)$ is coNP-complete. Section~\ref{sec:trubin} studies the Trubin property, and Sections~\ref{sec:diameter} and~\ref{sec:pivots} prove the diameter and pivot bounds. Appendices~\ref{app:adjacency}--\ref{app:pivots} contain the omitted proofs, in the same order.

\section{Adjacency of Multilinear Vertices}\label{sec:adjacency}
In this section we outline the proof of Theorem~\ref{thm:adjacency-sr} and show that adjacency in $\MP(H)$ is hard to decide. Let $S,T\subseteq V$ be distinct, and consider the midpoint $m=(z^S+z^T)/2$ of $z^S$ and $z^T$. Let $\mathcal W^m$ be the space of directions $w$ with $a^\top w=0$ for every inequality $a^\top z\le b$ of $\SR(H)$ tight at $m$. Equivalently, $w\in \mathcal{W}^m$ if and only if $m\pm\varepsilon w\in\SR(H)$ for some $\varepsilon>0$ (Lemma~\ref{lem:two-sided-directions} in the appendix). The space $\mathcal W^m$ is the translation space of the minimal face of $\SR(H)$ containing $m$, which is also the minimal face containing $z^S$ and $z^T$.

We denote by $\co(H)$ the number of connected components of a hypergraph $H$, and for a positive integer $k$ we write $[k]:=\{1,\dots,k\}$. Let $U:=S\triangle T$. For $v\notin U$, one of the bounds $0\le z_v\le1$ is tight at $m$, so every $w\in\mathcal W^m$ satisfies $w_v=0$. We call $S\setminus T$ and $T\setminus S$ the two \emph{sides} of $U$, and we define $\sigma(v):=1$ for $v\in S\setminus T$ and $\sigma(v):=-1$ for $v\in T\setminus S$. A generating edge $f$ forces $\sigma(u)w_u=\sigma(v)w_v$ for all $u,v\in f\cap U$. Indeed, in cases (a) and (b) of Definition~\ref{def:contracted-hypergraph}, we have $m_f=\frac12$, and the inequalities $z_f\le z_v$ with $v\in f\cap U$ are tight at $m$, so $w_u=w_v$. In case (c), we have $m_f=0$, and the inequalities $z_f\ge0$ and $z_f\ge\sum_{v\in f}z_v-|f|+1$ are tight at $m$, so $w_u+w_v=w_f=0$. Let $\pi$ denote the projection onto the coordinates indexed by $U$. Then $\pi(\mathcal W^m)$ consists exactly of the vectors $\lambda$ for which $\sigma(v)\lambda_v$ is constant on each connected component of $\Theta(H,S,T)$ (Lemma~\ref{lem:projection} in the appendix). Since every edge coordinate of a direction in $\mathcal W^m$ is zero or equal to one of its node coordinates indexed by $U$ (Lemma~\ref{lem:edge-direction} in the appendix), the projection $\pi$ is injective. Each component therefore contributes exactly one free parameter, which gives the following result.
\begin{proposition}\label{prop:minimal-face-dimension}
Let $H$ be a hypergraph and let $S,T\subseteq V$ be distinct. Then the minimal face of $\SR(H)$ containing $z^S$ and $z^T$ has dimension $\co(\Theta(H,S,T))$.
\end{proposition}
\noindent For graphs, Proposition~\ref{prop:minimal-face-dimension} is part of \cite[Theorem~8]{Pad89}. Since two distinct vertices are adjacent if and only if the minimal face containing them has dimension one, Proposition~\ref{prop:minimal-face-dimension} implies Theorem~\ref{thm:adjacency-sr}. For a set $A$, we denote by $\binom{A}{2}$ the collection of its subsets of cardinality two. By Theorem~\ref{thm:adjacency-sr}, the multilinear vertices are pairwise adjacent in $\SR(H)$ if and only if $\binom{V}{2}\subseteq E$. If $V\neq\emptyset$, this condition also characterizes $\diam(\MP(H))=1$ (Theorem~\ref{thm:complete-two-edge-layer} and Corollary~\ref{cor:diameter-one-mp} in the appendix).

Theorem~\ref{thm:adjacency-sr} describes adjacency in $\SR(H)$ combinatorially. In contrast, adjacency in $\MP(H)$ is hard to decide, even if every edge has cardinality three.
\begin{proposition}\label{prop:mp-hardness}
The problem of deciding, given a hypergraph $H$ and distinct $S,T\subseteq V$, whether $z^S$ and $z^T$ are adjacent in $\MP(H)$ is coNP-complete. It remains coNP-complete if every edge of $H$ has cardinality three.
\end{proposition}
\noindent\textit{Proof idea.} A convex combination of multilinear vertices that equals the midpoint of $z^S$ and $z^T$ and gives positive weight to a third vertex certifies non-adjacency, so the problem is in coNP. For hardness, the face of $\MP(H)$ on which all edge coordinates vanish is affinely isomorphic to the convex hull of the vectors $x\in\{0,1\}^V$ that satisfy $\sum_{v\in e}x_v\le|e|-1$ for every $e\in E$. Matsui proved that deciding non-adjacency of two given vertices of such a polytope is NP-complete when every edge has cardinality three \cite{Mat95} (see \cite[Corollary~4.2]{Mat94tr}).

\section{The Trubin Property}\label{sec:trubin}
In this section we study when $\SR(H)$ has the Trubin property with respect to $\MP(H)$. We begin with a smallest hypergraph for which the property fails.
\begin{proposition}\label{prop:sparse-counterexample}
Let $H=(V,E)$ with $V=\{1,2,3\}$ and $E=\{\{1,2\},\{1,2,3\}\}$, and let $S=\{1,2\}$ and $T=\{3\}$. Then $z^S$ and $z^T$ are adjacent in $\MP(H)$, but not in $\SR(H)$.
\end{proposition}
\noindent\textit{Proof idea.} In the coordinates $(z_1,z_2,z_3,z_{12},z_{123})$, consider the cost vector $c=(-1,-1,1,3,-2)$. Among the multilinear vertices, the only maximizers of $c^\top z$ are $z^S$ and $z^T$, so these two vertices are adjacent in $\MP(H)$. On the other hand, $\Theta(H,S,T)$ has the connected components $\{1,2\}$ and $\{3\}$, so by Theorem~\ref{thm:adjacency-sr} the vertices $z^S$ and $z^T$ are not adjacent in $\SR(H)$. Geometrically, their segment is a diagonal of a two-dimensional face of $\SR(H)$ rather than a one-dimensional face (Fig.~\ref{fig:counterexample}(b) and Remark~\ref{rem:counterexample-face} in the appendix).

Every failure of the Trubin property produces a fractional vertex, such as $\bar z$ in Fig.~\ref{fig:counterexample}(b) (Lemma~\ref{lem:fractional-failure} in the appendix). In contrast, the following condition on a pair of distinct sets $S,T\subseteq V$ keeps fractional vertices out of the minimal face of $\SR(H)$ containing $z^S$ and $z^T$:
\begin{itemize}
\item[$(\star)$] every edge $e\subseteq S\cup T$ that meets both $S\setminus T$ and $T\setminus S$ contains nodes $u\in S\setminus T$ and $v\in T\setminus S$ in the same connected component of $\Theta(H,S,T)$.
\end{itemize}
\begin{theorem}\label{thm:cube}
Let $H$ be a hypergraph, let $S,T\subseteq V$ be distinct and satisfy $(\star)$, and let $C_1,\dots,C_k$ be the connected components of $\Theta(H,S,T)$. Then the minimal face of $\SR(H)$ containing $z^S$ and $z^T$ is affinely isomorphic to the cube $[0,1]^k$, and its vertices are the $2^k$ multilinear vertices $z^{S\triangle C_I}$, where $C_I:=\bigcup_{l\in I}C_l$ for $I\subseteq[k]$. In particular, $z^S$ and $z^T$ are adjacent in $\MP(H)$ if and only if $\Theta(H,S,T)$ is connected.
\end{theorem}
\noindent\textit{Proof idea.} By the proof of Proposition~\ref{prop:minimal-face-dimension}, the points of the minimal face are parametrized around the midpoint of $z^S$ and $z^T$ by one coordinate $t_l$ for each component $C_l$. Condition $(\star)$ guarantees that all inequalities of $\SR(H)$ hold on the whole cube $t\in[-\frac12,\frac12]^k$, whose vertices correspond to the multilinear vertices $z^{S\triangle C_I}$.

When $H$ is a graph, every edge that meets both $S\setminus T$ and $T\setminus S$ is itself an edge of $\Theta(H,S,T)$, so condition $(\star)$ holds. Theorem~\ref{thm:cube} thus recovers Padberg's result that the minimal face of $\SR(G)$ containing $z^S$ and $z^T$ has only $0/1$ vertices \cite[Theorem~8]{Pad89}, and shows moreover that this face is a cube. In the counterexample of Proposition~\ref{prop:sparse-counterexample}, condition $(\star)$ fails for the edge $\{1,2,3\}$. Accordingly, the minimal face is not a cube, and it has the fractional vertex $\bar z$ (Fig.~\ref{fig:counterexample}(b)).

Let $S,T\subseteq V$ be distinct. Recall that the pair $S,T$ is bad if some edge $e\subseteq S\cup T$ has a trace that meets exactly two connected components of $\Theta(H,S,T)$. For example, in Proposition~\ref{prop:sparse-counterexample} the edge $\{1,2,3\}$ is its own trace and meets exactly the connected components $\{1,2\}$ and $\{3\}$, so the pair $S,T$ is bad. Every bad pair yields two multilinear vertices that violate the Trubin property.
\begin{theorem}\label{thm:bad-pair}
Let $H$ be a hypergraph, let $S,T\subseteq V$ be a bad pair, and let $e\subseteq S\cup T$ be an edge whose trace meets exactly the connected components $K$ and $K'$ of $\Theta(H,S,T)$. Then $z^{(S\cap T)\cup K}$ and $z^{(S\cap T)\cup K'}$ are adjacent in $\MP(H)$ but not in $\SR(H)$. In particular, if $H$ has a bad pair, then $\SR(H)$ does not have the Trubin property with respect to $\MP(H)$.
\end{theorem}
\noindent\textit{Proof idea.} Let $S':=(S\cap T)\cup K$ and $T':=(S\cap T)\cup K'$. The generating edges whose traces lie in $K$ or in $K'$ remain generating edges for the pair $S',T'$, and no generating edge for $S',T'$ meets both $K$ and $K'$. Hence $\Theta(H,S',T')$ has exactly the two connected components $K$ and $K'$. Now write the midpoint of $z^{S'}$ and $z^{T'}$ as a convex combination of multilinear vertices $z^R$. The generating edges force each $R$ with positive weight to contain all or none of $K$, and all or none of $K'$. The edge $e$ forbids $R\supseteq K\cup K'$. Hence the combination is unique, and $z^{S'}$ and $z^{T'}$ are adjacent in $\MP(H)$.

Let $\mathcal G(S,T)$ be the bipartite graph between the connected components of $\Theta(H,S,T)$ and the distinct traces that meet at least two of them, in which each such trace is joined to the components it meets. The following sufficient condition requires that the pair $S,T$ is not bad and that $\mathcal G(S,T)$ has no cycle.
\begin{theorem}\label{thm:acyclic-traces}
Let $H$ be a hypergraph and let $S,T\subseteq V$ be distinct. Suppose that the pair $S,T$ is not bad and that $\mathcal G(S,T)$ is a forest. Then $z^S$ and $z^T$ are adjacent in $\MP(H)$ if and only if $\Theta(H,S,T)$ is connected.
\end{theorem}
\noindent\textit{Proof idea.} Suppose that $\Theta(H,S,T)$ is disconnected. We show that some component $K$ can be switched in both $z^S$ and $z^T$ without changing their midpoint, so that $z^S$ and $z^T$ are not adjacent in $\MP(H)$. A trace prevents this switch only if its part in $K$ lies on one side of $S\triangle T$ and its remaining part lies on the other side. Since the pair $S,T$ is not bad, each trace prevents the switch of at most one component. Since $\mathcal G(S,T)$ is a forest, there are fewer such traces than components. Hence some component can be switched.

If the hypotheses of Theorem~\ref{thm:acyclic-traces} hold for all distinct $S,T\subseteq V$, then $\SR(H)$ has the Trubin property with respect to $\MP(H)$. By Theorem~\ref{thm:bad-pair}, the Trubin property fails as soon as the first hypothesis fails for some pair. Both hypotheses depend on the pair $S,T$. A condition on $H$ alone that implies them for every pair is the following. For every $e\in E$, let $M_e := \binom{e}{2}\setminus E$ be the set of pairs of nodes of $e$ that are not edges of $H$. We say that $H$ satisfies the \emph{matching condition} if $M_e$ is a matching for every $e\in E$.
\begin{theorem}\label{thm:matching-trubin}
    Let $H$ be a hypergraph that satisfies the matching condition. Then, for all distinct $S,T\subseteq V$, the vertices $z^S$ and $z^T$ are adjacent in $\MP(H)$ if and only if $\Theta(H,S,T)$ is connected. In particular, $\SR(H)$ has the Trubin property with respect to $\MP(H)$.
\end{theorem}
\noindent\textit{Proof idea.} We show that every trace lies in a single connected component of $\Theta(H,S,T)$. Let $u$ and $v$ be distinct nodes of the trace of an edge $e$. If $\{u,v\}\in E$, or if the trace is $\{u,v\}$, then $u$ and $v$ are joined by an edge of $\Theta(H,S,T)$. Otherwise, $\{u,v\}\in M_e$, and the trace has a node $r\notin\{u,v\}$. Since $M_e$ is a matching, $\{u,r\}$ and $\{v,r\}$ are edges of $H$, and hence of $\Theta(H,S,T)$. Hence no trace meets two or more components, and Theorem~\ref{thm:acyclic-traces} applies.

Since $M_e = \emptyset$ whenever $|e| = 2$, every graph satisfies the matching condition. Theorem~\ref{thm:matching-trubin} thus recovers the equivalence of (ii) and (iii) in Theorem~\ref{thm:connectivity-in-graph}, while allowing edges of arbitrary cardinality. Requiring the hypotheses of Theorem~\ref{thm:acyclic-traces} for all pairs is weaker than the matching condition. For example, one can check that the hypergraph with node set $\{1,2,3,4\}$ and edges $\{1,4\}$, $\{2,4\}$, $\{3,4\}$, $\{1,2,3\}$ satisfies them for every pair, although it violates the matching condition and is not Berge-acyclic. Condition $(\star)$ of Theorem~\ref{thm:cube} and the hypotheses of Theorem~\ref{thm:acyclic-traces} are incomparable (see the end of Appendix~\ref{app:trubin}).

For hypergraphs of rank at most three, bad pairs are the only obstruction to the Trubin property.
\begin{theorem}\label{thm:rank-three}
Let $H$ be a hypergraph of rank at most three. Then $\SR(H)$ has the Trubin property with respect to $\MP(H)$ if and only if $H$ has no bad pair.
\end{theorem}
\noindent\textit{Proof idea.} One direction is Theorem~\ref{thm:bad-pair}. For the other one, by Theorem~\ref{thm:acyclic-traces} it suffices to show that $\mathcal G(S,T)$ is a forest for every pair $S,T$. If the pair $S,T$ is not bad and $H$ has rank at most three, every trace that meets at least two components is an edge with one node in each of three components. Consider a shortest cycle of $\mathcal G(S,T)$. If it consists of two traces that meet the same three components, adding a suitable node to both $S$ and $T$ yields a bad pair. Otherwise, putting one component of the cycle alone on one side of $S\triangle T$ yields a bad pair.

By Theorem~\ref{thm:rank-three}, for hypergraphs of rank at most three, a failure of the Trubin property is certified by a pair $S,T$ together with an edge whose trace meets exactly two connected components of $\Theta(H,S,T)$. This certificate can be checked in polynomial time. Hence, for these hypergraphs, deciding whether $\SR(H)$ has the Trubin property with respect to $\MP(H)$ is in coNP.

\section{Diameter Bounds}\label{sec:diameter}
In this section we bound the diameters of $\MP(H)$ and $\SR(H)$. For two vertices $p$ and $q$ of a polytope $P$, the \emph{distance} $\dist_P(p,q)$ is the minimum number of edges of a path between them in the one-skeleton of $P$. Let $S,T\subseteq V$ be distinct, and let $C_1,\ldots,C_k$ be the connected components of $\Theta(H,S,T)$ in an arbitrary order. Set $S_0:=S$ and $S_i:=S_{i-1}\triangle C_i$ for $i\in[k]$, so that $S_k=T$. Every generating edge $f$ with $f\cap(S\triangle T)\subseteq C_i$ is also a generating edge for the pair $S_{i-1},S_i$, with the same trace. These traces connect $C_i$, which is the node set of $\Theta(H,S_{i-1},S_i)$. Hence $\Theta(H,S_{i-1},S_i)$ is connected. By Theorem~\ref{thm:adjacency-sr}, the vertices $z^{S_{i-1}}$ and $z^{S_i}$ are adjacent in $\SR(H)$, and hence in $\MP(H)$. Since the components are nonempty and pairwise disjoint subsets of $S\triangle T$, this yields the following bound.
\begin{theorem}\label{thm:diameter-mp}
    Let $H$ be a hypergraph and let $S,T\subseteq V$ be distinct. Then
    \[
        \dist_{\SR(H)}(z^S,z^T)\le \co(\Theta(H,S,T)) \le |V|,
    \]
    and the same bounds hold for $\dist_{\MP(H)}(z^S,z^T)$. In particular, $\diam(\MP(H))\le |V|$, and this bound is attained when $E=\emptyset$.
\end{theorem}
\noindent The bound $\diam(\MP(H))\le|V|$ improves the bound $|V|+|E|$ given by Naddef's theorem that the diameter of a $0/1$ polytope is at most its dimension \cite{Nad89}.

Theorem~\ref{thm:diameter-mp} concerns only multilinear vertices. For graphs, Padberg's half-integrality theorem \cite[Theorem~7 and p.~157]{Pad89} allows a fractional vertex to be rounded one connected component of its fractional nodes at a time. This gives the same bound as Theorem~\ref{thm:diameter-mp}.
\begin{theorem}\label{thm:diameter-sr-graph}
Let $G=(V,E)$ be a graph. Then the diameter of $\SR(G)$ is at most $|V|$, and this bound is attained when $E=\emptyset$.
\end{theorem}
\noindent\textit{Proof idea.} We proceed by induction on $|V|$. Fixing a node coordinate to $0$ or $1$ gives a face of $\SR(G)$ that is affinely isomorphic to the standard relaxation of a graph with one fewer node. Hence two vertices that agree on an integral node coordinate are at distance at most $|V|-1$. If both vertices are multilinear, Theorem~\ref{thm:diameter-mp} applies. Otherwise, we round a connected component of the fractional nodes of one of the two vertices, or of both, along a spanning tree. Each rounding is a single step, and the roundings can be chosen so that the two vertices come to agree on at least as many integral node coordinates as the number of steps used.

Already for hypergraphs of rank three, the vertices of $\SR(H)$ need not be half-integral, and their denominators can even be exponentially large.
\begin{proposition}\label{prop:denominators}
For every integer $\ell\ge2$, there is a hypergraph $H$ of rank three with $|V|=|E|=2\ell+1$ such that $\SR(H)$ has a vertex with a coordinate whose denominator, in lowest terms, is $2^\ell+(-1)^\ell$.
\end{proposition}
\noindent\textit{Proof idea.} The hypergraph is a cycle on $2\ell+1$ nodes, with the edge $\{1,2\}$, the pairs $\{2r+1,2r+2\}$ and the triples $\{2r,2r+1,2r+2\}$ for $r\in[\ell]$, where node $2\ell+2$ is node $1$. At the vertex, the tight constraints force $z_{2r+1}=z_{2r+2}$ and $z_{2r}+z_{2r+1}+z_{2r+2}=2$. Hence $u_r:=z_{2r}$ satisfies $u_{r+1}=1-u_r/2$ along the cycle. Closing the cycle with $z_1+z_2=1$ gives a unique solution, whose denominator is $2^\ell+(-1)^\ell$.

For a positive integer $D$, if all vertices of an $n$-dimensional polytope have coordinates in $\{0,\frac1D,\frac2D,\dots,1\}$, its diameter is at most $Dn$ \cite{KleOnn92}, and at most $(D-\frac12)n$ for $D\ge2$ \cite{dPMic16DCG}. By Proposition~\ref{prop:denominators}, these bounds can be exponential in $|V|$ for $\SR(H)$.

For arbitrary hypergraphs, we work instead with bases of tight constraints. We call \emph{rows} the equations obtained by setting the defining inequalities of $\SR(H)$ to equality, and we denote them by
\[
    \begin{array}{lll}
    N_v^0: z_v=0, \qquad & N_v^1: z_v=1, \qquad & A^e: z_e=0,\\[2pt]
    B_v^e: z_e=z_v \ \text{ for } v\in e, \qquad & C^e: z_e=\sum_{v\in e}z_v-|e|+1. &
    \end{array}
\]
We call $N_v^0$ and $N_v^1$ the \emph{node bounds} of $v$. A row is \emph{tight} at a point of $\SR(H)$ if the point satisfies it. Let $d:=|V|+|E|$. A \emph{basis} at a vertex $z$ of $\SR(H)$ is a set of $d$ linearly independent rows tight at $z$. A \emph{pivot} replaces one row of a basis $\mathcal B$ at a vertex $z$ by another row so that the result is a basis $\mathcal B'$ at a vertex $z'$. The pivot is \emph{degenerate} if $z'=z$, and \emph{nondegenerate} otherwise. For two bases $\mathcal B$ and $\mathcal B'$, we denote by $\piv(\mathcal B,\mathcal B')$ the minimum number of pivots that transform $\mathcal B$ into $\mathcal B'$. Since pivots are reversible, $\piv(\mathcal B,\mathcal B')=\piv(\mathcal B',\mathcal B)$. In a nondegenerate pivot, the $d-1$ rows common to $\mathcal B$ and $\mathcal B'$ are linearly independent and tight at $z$ and $z'$, so the face of $\SR(H)$ on which they are tight is a segment of the line they define. Since $z$ and $z'$ are distinct vertices of $\SR(H)$ in this face, they are its endpoints, and the face is the edge $[z,z']$ of the one-skeleton. Hence a sequence of $t$ pivots moves along at most $t$ edges of the one-skeleton of $\SR(H)$. By the basis exchange property of linear matroids, a basis $\mathcal B$ can be transformed into any other basis $\mathcal B'$ at the same vertex by $|\mathcal B'\setminus\mathcal B|\le d$ degenerate pivots (Lemma~\ref{lem:basis-exchange} in the appendix).

Let $\mathcal B$ be a basis. We call an edge $e$ \emph{low} if the rows of $\mathcal B$ involving $z_e$ are exactly $A^e$ and $C^e$, that is, the two lower bounds on $z_e$. A basis without low edges is at a multilinear vertex (Lemma~\ref{lem:active-system-integrality} in the appendix). Indeed, its rows $B^e_v$ force equal values on groups of node coordinates, and the remaining rows fix each group to $0$ or $1$. Since a basis has $|V|+|E|$ rows and every edge coordinate is involved in at least one of them, a basis has at most $\min\{|V|,|E|\}$ low edges. To remove a low edge $e$, we drop $C^e$ and move along the line defined by the other rows until a new inequality becomes tight. With at most one preliminary pivot for every edge that is not low, each low edge is removed with at most two pivots. This yields the following reduction.
\begin{lemma}\label{lem:frac-int}
For every vertex $z$ of $\SR(H)$ and every basis $\mathcal B$ at $z$, there exist a multilinear vertex $\hat z$ and a basis $\hat{\mathcal B}$ at $\hat z$ such that
\[
\piv(\mathcal B,\hat{\mathcal B})\le|E|+\min\{|V|,|E|\},
\qquad \dist_{\SR(H)}(z,\hat z)\le|E|+\min\{|V|,|E|\}.
\]
\end{lemma}
\noindent Applying Lemma~\ref{lem:frac-int} at both endpoints, and Theorem~\ref{thm:diameter-mp} between the two resulting multilinear vertices, gives the following bound.
\begin{theorem}\label{thm:diameter-sr}
    Let $H$ be a hypergraph. Then
    \[
        \diam(\SR(H))\le |V|+2|E|+2\min\{|V|,|E|\}.
    \]
\end{theorem}

\section{Pivot Bounds}\label{sec:pivots}
In this section we bound the number of pivots between bases at two vertices of $\SR(H)$. We first consider two multilinear vertices at which the bases may be chosen. The following bound does not depend on $|E|$. Its proof switches one node at a time instead of whole components.
\begin{theorem}\label{thm:pivot-flip}
Let $H$ be a hypergraph and let $S,T\subseteq V$. Then there exist bases $\mathcal B^S$ at $z^S$ and $\mathcal B^T$ at $z^T$, and a sequence of at most $|S\triangle T|$ pivots from $\mathcal B^S$ to $\mathcal B^T$ that visits only bases at multilinear vertices. In particular,
\[
\piv(\mathcal B^S,\mathcal B^T)\le|S\triangle T|\le|V|,
\]
and this bound is attained when $E=\emptyset$. Moreover, for every basis $\mathcal B$ at $z^S$, there is a basis $\mathcal B^T$ at $z^T$ with $\piv(\mathcal B,\mathcal B^T)\le d+|S\triangle T|\le 2|V|+|E|$.
\end{theorem}
\noindent\textit{Proof idea.} We remove the nodes of $S\setminus T$ one at a time, and then add the nodes of $T\setminus S$ one at a time. For each edge $e$, one of the rows $A^e$ and $B^e_w$ with $w\in e$ is tight along this whole path of multilinear vertices. For example, if $e\subseteq S$ and $e\not\subseteq T$, then $B^e_u$ is tight throughout, where $u$ is the first removed node of $e$. Together with one tight node bound for each node, these rows form a basis, and consecutive bases differ only in the node bound of the switched node. The second bound follows from Lemma~\ref{lem:basis-exchange}.

Switching whole components instead of single nodes gives bases at $z^S$ and $z^T$ that are at most $\co(\Theta(H,S,T))$ pivots apart whenever the pair $S,T$ satisfies condition $(\star)$ (Theorem~\ref{thm:pivot-integral} in the appendix).

If the initial basis is prescribed, the dependence on $|E|$ cannot be avoided, even between adjacent vertices.
\begin{proposition}\label{prop:lowerbound}
For every integer $n\ge1$, there is a hypergraph $H=(V,E)$ with $|V|=n+1$ and $|E|=2^n-1$, two multilinear vertices $z^S$ and $z^T$ that are adjacent in $\SR(H)$, and a basis $\mathcal B^S$ at $z^S$ such that $\piv(\mathcal B^S,\mathcal B^T)\ge|E|+1$ for every basis $\mathcal B^T$ at $z^T$. Equality is attained for a suitable $\mathcal B^T$.
\end{proposition}
\noindent\textit{Proof idea.} Let $V=\{a,b_1,\ldots,b_n\}$, $E=\{e\subseteq V:a\in e,\ |e|\ge2\}$, $S=V$, and $T=V\setminus\{a\}$. Let $\mathcal B^S$ consist of the rows $N_v^1$ for $v\in V$ and, for every $e\in E$, a row $B^e_{b(e)}$ with $b(e)\in e\setminus\{a\}$. None of the $|E|+1$ rows $B^e_{b(e)}$ and $N_a^1$ is tight at $z^T$, and each pivot exchanges one row.

Finally, we bound the number of pivots between prescribed bases at arbitrary vertices of $\SR(H)$.
\begin{theorem}\label{thm:pivot-sr}
    Let $H$ be a hypergraph, let $p$ and $q$ be vertices of $\SR(H)$, and let $\mathcal{B}$ and $\mathcal{B}'$ be bases at $p$ and $q$, respectively. Then
    \[
        \piv(\mathcal{B},\mathcal{B}') \le 3|V| + 4|E|+2\min\{|V|,|E|\}.
    \]
\end{theorem}
\noindent\textit{Proof idea.} By Lemma~\ref{lem:frac-int}, each prescribed basis reaches a basis at a multilinear vertex in at most $|E|+\min\{|V|,|E|\}$ pivots. Theorem~\ref{thm:pivot-flip} connects these two vertices, with at most $d$ degenerate pivots at each junction by Lemma~\ref{lem:basis-exchange}. In total, this gives at most $2(|E|+\min\{|V|,|E|\})+2d+|V|=3|V|+4|E|+2\min\{|V|,|E|\}$ pivots.

\bigskip
\begin{small}\sloppy
\paragraph{Funding.} A. Del Pia is partially supported by AFOSR Grant \mbox{FA9550-23-1-0433} and ONR Grant \mbox{N00014-25-1-2490}. Any opinions, findings, conclusions, or recommendations expressed in this material are those of the authors and do not necessarily reflect the views of the Air Force Office of Scientific Research or the Office of Naval Research.

\paragraph{AI methodology.} AI tools (Claude, developed by Anthropic) were used for proofreading and suggesting improvements to the exposition. The authors take full responsibility for the content of the paper.
\end{small}

{\small
\bibliographystyle{plainurl}
\bibliography{biblio}
}
\clearpage
\appendix
\section{Adjacency of Multilinear Vertices}\label{app:adjacency}
This appendix contains the proofs omitted from the main text. Appendices~\ref{app:adjacency}--\ref{app:pivots} correspond to Sections~\ref{sec:adjacency}--\ref{sec:pivots}. We begin with some polyhedral preliminaries. Let $P\subseteq\mathbb{R}^N$ be a polytope described by
\[
    P=\{x\in\mathbb R^N:a_h^\top x\le b_h,\ h\in \mathcal I\}.
\]
For $p\in P$, denote by
\[
    \mathcal T(p):=\{h\in \mathcal I:a_h^\top p=b_h\}
\]
the index set of inequalities tight at $p$. It is standard that the unique minimal face of $P$ containing $p$ is
\[
    F^p := \{x\in P:a_h^\top x = b_h,\ \forall h\in \mathcal T(p)\},
\]
and that its affine hull is $p+\mathcal W^p$, where
\[
    \mathcal W^p :=
    \{w\in\mathbb R^N:a_h^\top w=0,\ \forall h\in\mathcal T(p)\}.
\]
Indeed, $F^p$ is a face of $P$, as the intersection of the faces induced by the inequalities in $\mathcal T(p)$, and it contains $p$. Since every inequality not in $\mathcal T(p)$ has positive slack at $p$ and $\mathcal I$ is finite, a sufficiently small neighborhood of $p$ in the affine space $p+\mathcal W^p$ is contained in $F^p$. Therefore $p+\mathcal W^p$ is the affine hull of $F^p$, and $p$ lies in the relative interior of $F^p$. Hence every face of $P$ containing $p$ contains $F^p$.

\begin{lemma}\label{lem:two-sided-directions}
    Given $p\in P$, the translation space $\mathcal{W}^p$ is the set of all two-sided feasible directions at $p$ in $P$. That is, 
    \[
        \mathcal W^p = \{w\in\mathbb{R}^N:\exists \varepsilon>0\text{ such that }p\pm \varepsilon w\in P\}.
    \]
\end{lemma}

\begin{proof}
Let $w\in\mathcal W^p$. For every $h\in\mathcal T(p)$, we have $a_h^\top w=0$, and hence $a_h^\top(p\pm\varepsilon w)=b_h$ for every $\varepsilon>0$. For every $h\notin\mathcal T(p)$, the slack $b_h-a_h^\top p$ is positive, so $a_h^\top(p\pm\varepsilon w)\le b_h$ for sufficiently small $\varepsilon>0$. Thus $p\pm\varepsilon w\in P$ for sufficiently small $\varepsilon>0$. Conversely, let $w\in\mathbb R^N$ and $\varepsilon>0$ satisfy $p\pm\varepsilon w\in P$, and let $h\in\mathcal T(p)$. Since $a_h^\top p=b_h$, the feasibility of $p+\varepsilon w$ and of $p-\varepsilon w$ gives $a_h^\top w\le0$ and $a_h^\top w\ge0$. Hence $a_h^\top w=0$ for every $h\in\mathcal T(p)$, that is, $w\in\mathcal W^p$.
\end{proof}

\begin{lemma}\label{lem:adjacency-subpolytope}
    Let $P\subseteq Q$ be polytopes. If two vertices $p,q$ of $P$ are adjacent in $Q$, then they are adjacent in $P$.
\end{lemma}

\begin{proof}
    Since $p$ and $q$ are adjacent in $Q$, there exists a linear function $c$ such that
    \[
        \conv\{p,q\} = \{x\in Q:c^\top x=\max_{y\in Q}c^\top y\}.
    \]
    Since $P\subseteq Q$ and $p,q\in P$, the maximum of $c^\top x$ over $P$ is the same, and
    \[
        P\cap\conv\{p,q\} = \{x\in P:c^\top x=\max_{y\in P}c^\top y\}.
    \]
    Because $p,q\in P$ and $P$ is convex, $\conv\{p,q\}\subseteq P$. Hence
    \[
        P\cap\conv\{p,q\}=\conv\{p,q\}.
    \]
    Therefore $\conv\{p,q\}$ is a one-dimensional face of $P$, and thus $p$ and $q$ are adjacent in $P$.
\end{proof}

The following characterization of non-adjacency is a standard fact. We include a proof for completeness.
\begin{lemma}\label{lem:midpoint-certificate}
Let $P$ be a polytope, let $p^1$ and $p^2$ be distinct vertices of $P$, and let $m:=\frac12(p^1+p^2)$. Then $p^1$ and $p^2$ are not adjacent in $P$ if and only if $m$ is a convex combination of vertices of $P$ in which some vertex other than $p^1$ and $p^2$ has a positive coefficient.
\end{lemma}

\begin{proof}
Suppose that $p^1$ and $p^2$ are not adjacent, and let $F$ be the minimal face of $P$ containing them. Since $F$ is not the segment between $p^1$ and $p^2$, it has a vertex different from $p^1$ and $p^2$. Since $m=\frac12p^1+\frac12p^2$, a face of $P$ contains $m$ if and only if it contains $p^1$ and $p^2$. Hence $F$ is the minimal face of $P$ containing $m$, and $m$ lies in the relative interior of $F$ by the preliminaries above. Let $b$ be the average of the vertices of $F$. Then $m':=m+\varepsilon(m-b)\in F$ for some $\varepsilon>0$, and $m=\frac{1}{1+\varepsilon}m'+\frac{\varepsilon}{1+\varepsilon}b$. Writing $m'$ as a convex combination of the vertices of $F$, we obtain a convex combination equal to $m$ in which every vertex of $F$ has a positive coefficient.

Conversely, suppose that such a combination exists and that $p^1$ and $p^2$ are adjacent. Then there is a vector $c$ such that the maximizers of $c^\top x$ over $P$ are exactly the points of the segment between $p^1$ and $p^2$. Since $m$ is a maximizer, every vertex with a positive coefficient in the combination is a maximizer. Hence it lies in the segment, and as a vertex of $P$ it is $p^1$ or $p^2$, a contradiction.
\end{proof}

\begin{lemma}\label{lem:midpoint-nonadjacency}
    Let $P$ be a polytope, and let $p^1,p^2,q^1,q^2$ be vertices of $P$ such that
    \[
        \{p^1,p^2\}\neq\{q^1,q^2\}
    \]
    and
    \[
        p^1+p^2=q^1+q^2.
    \]
    Then $p^1$ and $p^2$ are not adjacent in $P$.
\end{lemma}

\begin{proof}
If $q^i=p^j$ for some $i,j\in\{1,2\}$, then $p^1+p^2=q^1+q^2$ shows that the remaining two vertices are also equal, so $\{q^1,q^2\}=\{p^1,p^2\}$, a contradiction. Hence $q^1$ and $q^2$ are different from $p^1$ and $p^2$. If $p^1=p^2$, there is nothing to prove. Otherwise, $\frac12(p^1+p^2)=\frac12q^1+\frac12q^2$ is a convex combination of vertices of $P$ in which a vertex other than $p^1$ and $p^2$ has a positive coefficient. By Lemma~\ref{lem:midpoint-certificate}, $p^1$ and $p^2$ are not adjacent in $P$.
\end{proof}

\subsection{Tight constraints and the contracted hypergraph}\label{sec:app-adjacency}
Throughout this subsection, let $S,T\subseteq V$ be distinct, let $m=(z^S+z^T)/2$ and
$U=S\triangle T$, and let $\mathcal W^m$ be the space of Lemma~\ref{lem:two-sided-directions}
with $p=m$. Let $\sigma(v):=1$ for $v\in S\setminus T$ and $\sigma(v):=-1$ for $v\in T\setminus S$. For every $v\notin U$ we have $m_v\in\{0,1\}$, so a node bound of $v$ is tight
at $m$. Hence $w_v=0$ for every $w\in\mathcal W^m$ and every $v\notin U$.

\begin{lemma}\label{lem:active-edge-coupling}
Let $w\in\mathcal W^m$ and $f\in F_\Theta$. Then $\sigma(u)w_u=\sigma(v)w_v$ for all
$u,v\in f\cap U$.
\end{lemma}

\begin{proof}
Suppose first that $f$ satisfies condition (a) of Definition~\ref{def:contracted-hypergraph}.
Then $f\subseteq S$ and $f\not\subseteq T$, so $m_f=\frac12$. Every $v\in f\cap U$ lies in
$S\setminus T$, so $m_v=\frac12$ and the inequality $z_f\le z_v$ is tight at $m$. Hence
$w_f=w_v$ for every $v\in f\cap U$, and all these coordinates are equal. Since $\sigma$ is
constant on $f\cap U$, the claim follows. Condition (b) is symmetric.
Suppose now that $f$ satisfies condition (c), and let $f\cap U=\{u,v\}$ with
$u\in S\setminus T$ and $v\in T\setminus S$. Then $f\setminus U\subseteq S\cap T$, and $f$ is
contained in neither $S$ nor $T$, so $m_f=0$. The inequality $z_f\ge0$ is tight at $m$, which
gives $w_f=0$. Moreover, $\sum_{x\in f}m_x-|f|+1=(|f|-2)+1-|f|+1=0=m_f$, so the inequality
$z_f\ge\sum_{x\in f}z_x-|f|+1$ is tight at $m$ as well. Since $w_x=0$ for $x\in f\setminus U$,
this gives $0=w_f=w_u+w_v$, that is, $\sigma(u)w_u=\sigma(v)w_v$.
\end{proof}

\begin{lemma}\label{lem:edge-direction}
For every $w\in \mathcal{W}^m$ and every edge $f\in E$, at least one of the following holds: $w_f=0$ or there exists a node $v\in f\cap U$ such that $w_f=w_v$. 
\end{lemma}

\begin{proof}
    Fix $f\in E$. Since $z_f^S,z_f^T\in \{0,1\}$, we have
    \[
        m_f=\frac{1}{2}(z_f^S+z_f^T)\in\{0,\frac{1}{2},1\}.
    \]
    First suppose $m_f=0$. Then the lower bound $z_f\ge 0$ is tight at $m$. Therefore every direction $w\in \mathcal{W}^m$ satisfies $w_f=0$. If $m_f=1$, then $z_f^S=z_f^T=1$, so $f\subseteq S\cap T$. For every $v\in f$, we have $m_f=m_v=1$, so the inequality $z_f\le z_v$ is tight, and passing to directions gives $w_f=w_v=0$. Now suppose $m_f=\frac{1}{2}$. Then exactly one of $z_f^S$ and $z_f^T$ is equal to one. Equivalently, $f$ is contained in exactly one of $S$ and $T$. If $f\subseteq S$ and $f\not\subseteq T$, then there exists $v\in f\cap(S\setminus T)$. For this node, $m_f= m_v=\frac{1}{2}$ implies $w_f=w_v$. The case $f\subseteq T$ and $f\not\subseteq S$ is symmetric and gives $w_f=w_v$ for some $v\in f\cap(T\setminus S)$. Thus, in every case, $w_f$ is either zero or equal to a node coordinate indexed by $U$.
\end{proof}

\begin{lemma}\label{lem:projection}
    Let $C_1,\ldots,C_k$ be the connected components of $\Theta(H,S,T)$. Define 
    \[
        \pi:\mathcal{W}^m\longrightarrow \mathbb{R}^{U},\quad
        \pi(w)=(w_v)_{v\in U}.
    \]
    Then
    \begin{equation}\tag{$\dagger$}\label{eq:projection-image}
        \pi(\mathcal{W}^m)=\{\lambda \in \mathbb{R}^U:\sigma(u)\lambda_u = \sigma(v)\lambda_v \, \forall u,v\in C_l, l\in [k]\}.
    \end{equation}
\end{lemma}

\begin{proof}
Let $\mathcal R$ denote the right-hand side of \eqref{eq:projection-image}. We first show that $\pi(\mathcal W^m)\subseteq\mathcal R$. Let $w\in\mathcal W^m$. By Lemma~\ref{lem:active-edge-coupling}, we have $\sigma(u)w_u=\sigma(v)w_v$ for every edge $e\in E_\Theta$ and all $u,v\in e$. Hence $\sigma(v)w_v$ is constant along every path in $\Theta(H,S,T)$, and therefore on each connected component $C_l$. Thus $\pi(\mathcal W^m)\subseteq\mathcal R$.

Conversely, let $\lambda\in\mathcal R$. Define $w$ by $w_v=\lambda_v$ for $v\in U$ and $w_v=0$ for $v\notin U$. For each $f\in E$, set $w_f=0$ if $m_f\in\{0,1\}$, and $w_f=w_v$ for an arbitrary $v\in f\cap U$ if $m_f=\frac12$. This is well defined. Indeed, if $m_f=\frac12$, then $f$ is contained in exactly one of $S$ and $T$, so $f\cap U$ is nonempty and $\sigma$ is constant on $f\cap U$. If moreover $|f\cap U|\ge2$, then $f$ satisfies condition (a) or (b) of Definition~\ref{def:contracted-hypergraph}, so $f\cap U\in E_\Theta$ and all nodes of $f\cap U$ lie in the same connected component. Since $\lambda\in\mathcal R$, their $\lambda$-coordinates are equal.

It remains to show that $w\in\mathcal W^m$, that is, that $w$ satisfies the homogeneous equation of every inequality tight at $m$. If a node bound of $v$ is tight, then $v\notin U$ and $w_v=0$. If $z_f\ge0$ is tight, then $m_f=0$ and $w_f=0$. Suppose that $z_f\le z_v$ is tight. If $m_f=m_v\in\{0,1\}$, then $v\notin U$ and $w_f=0=w_v$. If $m_f=m_v=\frac12$, then $v\in f\cap U$ and $w_f=w_v$ by construction. Finally, suppose that $z_f\ge\sum_{v\in f}z_v-|f|+1$ is tight, so that we need $w_f=\sum_{v\in f}w_v$. If $m_f=1$, then $f\subseteq S\cap T$, and both sides are zero. If $m_f=\frac12$, the inequality is tight only when $f\cap U=\{u\}$ for some node $u$ and $f\setminus U\subseteq S\cap T$, and then $w_f=w_u=\sum_{v\in f}w_v$. If $m_f=0$, the inequality is tight in two cases only. In the first case, $f\cap U=\emptyset$ and exactly one node of $f$ lies outside $S\cup T$, and both sides are zero. In the second case, $f\subseteq S\cup T$ and $f\cap U=\{u,v\}$ with $u\in S\setminus T$ and $v\in T\setminus S$. Then $f$ satisfies condition (c), so $\{u,v\}\in E_\Theta$. Since $\lambda\in\mathcal R$, we have $\lambda_u=\sigma(u)\lambda_u=\sigma(v)\lambda_v=-\lambda_v$, and hence $\sum_{x\in f}w_x=w_u+w_v=0=w_f$. Thus $w\in\mathcal W^m$ and $\pi(w)=\lambda$, which proves $\mathcal R\subseteq\pi(\mathcal W^m)$.
\end{proof}

\begin{proof}[Proof of Proposition~\ref{prop:minimal-face-dimension}]
Since $m=\frac12 z^S+\frac12 z^T$, a face of $\SR(H)$ contains $m$ if and only if it contains both $z^S$ and $z^T$. Hence the minimal face $F^m$ of $\SR(H)$ containing $m$ is also the minimal face containing $z^S$ and $z^T$, and it suffices to show that $\dim(F^m)=\co(\Theta(H,S,T))$. Consider the linear map $\pi$ of Lemma~\ref{lem:projection}, and let $w\in\mathcal W^m$ satisfy $\pi(w)=0$. Then $w_v=0$ for every $v\in U$, and also for every $v\notin U$. By Lemma~\ref{lem:edge-direction}, every edge coordinate of $w$ is zero or equal to a node coordinate, so $w=0$. Hence $\pi$ is injective, and $\dim(\mathcal W^m)=\dim(\pi(\mathcal W^m))$. Let $k:=\co(\Theta(H,S,T))$. By Lemma~\ref{lem:projection}, the vectors in $\pi(\mathcal W^m)$ are exactly the vectors $\lambda$ with $\lambda_v=\sigma(v)t_l$ for all $v\in C_l$ and $l\in[k]$, where $(t_1,\dots,t_k)\in\mathbb R^k$. Thus $\dim(\pi(\mathcal W^m))=k$. Since $\mathcal W^m$ is the translation space of the affine hull of $F^m$, we conclude that $\dim(F^m)=\dim(\mathcal W^m)=k$.
\end{proof}

\begin{proof}[Proof of Theorem~\ref{thm:adjacency-sr}]
Let $m=\frac{1}{2}(z^S+z^T)$, and let $F^m$ be the minimal face of $\SR(H)$ containing $z^S$ and $z^T$. The vertices $z^S$ and $z^T$ are adjacent in $\SR(H)$ if and only if $\dim(F^m)=1$. By Proposition~\ref{prop:minimal-face-dimension}, $\dim(F^m)=\co(\Theta(H,S,T))$. Hence $z^S$ and $z^T$ are adjacent in $\SR(H)$ if and only if $\co(\Theta(H,S,T))=1$, which is exactly the statement that $\Theta(H,S,T)$ is connected.
\end{proof}

Theorem~\ref{thm:adjacency-sr} determines when all multilinear vertices are pairwise adjacent.
\begin{theorem}\label{thm:complete-two-edge-layer}
Let $H$ be a hypergraph. Then every pair of distinct multilinear vertices is adjacent in $\SR(H)$ if and only if $\binom{V}{2}\subseteq E$.
\end{theorem}
\begin{proof}
Suppose first that $\binom{V}{2}\subseteq E$, and let $S,T\subseteq V$ be distinct. If $|S\triangle T|=1$, then $\Theta(H,S,T)$ has a single node and is connected. Otherwise, let $u,v\in S\triangle T$ be distinct. Since $\binom{V}{2}\subseteq E$, we have $\{u,v\}\in E$. If $u$ and $v$ lie on the same side of the symmetric difference, then $\{u,v\}$ satisfies condition (a) or (b) of Definition~\ref{def:contracted-hypergraph}. Otherwise, it satisfies condition (c). Hence $\{u,v\}\in E_\Theta$. Thus $\Theta(H,S,T)$ is connected, and $z^S$ and $z^T$ are adjacent in $\SR(H)$ by Theorem~\ref{thm:adjacency-sr}.

Conversely, suppose that $\{u,v\}\notin E$ for some distinct $u,v\in V$, and let $S:=\{u,v\}$ and $T:=\emptyset$. Then $S\triangle T=S\cup T=\{u,v\}$. Since every edge of $H$ has cardinality at least two and $\{u,v\}\notin E$, no edge of $H$ is contained in $S\cup T$. Hence $E_\Theta=\emptyset$, and $\Theta(H,S,T)$ consists of the two isolated nodes $u$ and $v$. By Theorem~\ref{thm:adjacency-sr}, $z^S$ and $z^T$ are not adjacent in $\SR(H)$.
\end{proof}

\begin{proof}[Proof of Proposition~\ref{prop:mp-hardness}]
Let $m:=\frac12(z^S+z^T)$. We first show that the problem is in coNP. Consider the linear program
\[
\max\Bigl\{\sum_{R\subseteq V,\ R\neq S,T}\lambda_R\ :\ \sum_{R\subseteq V}\lambda_Rz^R=m,\ \ \sum_{R\subseteq V}\lambda_R=1,\ \ \lambda\ge0\Bigr\}
\]
in the variables $\lambda_R$ for $R\subseteq V$. It is feasible and bounded, and it has $d+1$ equality constraints. Hence it has an optimal basic solution, which has at most $d+1$ positive coefficients. These coefficients solve a nonsingular linear system with $0/1$ coefficients and a half-integral right-hand side, so by Cramer's rule their encoding size is polynomial in $d$. If $z^S$ and $z^T$ are not adjacent, the optimum is positive by Lemma~\ref{lem:midpoint-certificate}, and such a basic solution is a certificate of non-adjacency that can be checked in polynomial time.

For hardness, let $F_0$ be the face of $\MP(H)$ on which $z_e=0$ for every $e\in E$. Its vertices are the multilinear vertices $z^R$ such that no edge of $H$ is contained in $R$, that is, such that $\sum_{v\in e}z^R_v\le|e|-1$ for every $e\in E$. Since all edge coordinates vanish on $F_0$, the projection onto the node coordinates maps $F_0$ by an affine bijection onto
\[
P_H:=\conv\Bigl\{x\in\{0,1\}^V:\textstyle\sum_{v\in e}x_v\le|e|-1\ \ \forall e\in E\Bigr\}.
\]
A face of $F_0$ is a face of $\MP(H)$, and a face of $\MP(H)$ that contains two vertices of $F_0$ and has dimension one is contained in $F_0$. Hence two vertices of $F_0$ are adjacent in $\MP(H)$ if and only if their projections are adjacent in $P_H$. Let $\mathbf 1$ denote the all-ones vector. Matsui \cite{Mat95} (see \cite[Corollary~4.2]{Mat94tr}) proved that the following problem is NP-complete. Given a $0/1$ matrix $A$ with $n$ columns whose rows have exactly three ones each, and two distinct vectors $x^1,x^2\in\{0,1\}^n$ with $Ax^i\le 2\cdot\mathbf 1$ for $i=1,2$, decide whether $x^1$ and $x^2$ are non-adjacent vertices of $\conv\{x\in\{0,1\}^n:Ax\le2\cdot\mathbf 1\}$. Given such an instance, let $H$ be the hypergraph with node set $[n]$ whose edges are the supports of the distinct rows of $A$, and let $S$ and $T$ be the supports of $x^1$ and $x^2$. Then every edge of $H$ has cardinality three, $P_H=\conv\{x\in\{0,1\}^n:Ax\le2\cdot\mathbf 1\}$, and $z^S$ and $z^T$ are the vertices of $F_0$ whose projections are $x^1$ and $x^2$. Hence $x^1$ and $x^2$ are non-adjacent in $P_H$ if and only if $z^S$ and $z^T$ are not adjacent in $\MP(H)$. Since $H$, $S$, and $T$ are computed in polynomial time, deciding adjacency in $\MP(H)$ is coNP-hard, even if every edge of $H$ has cardinality three.
\end{proof}

\section{The Trubin Property}\label{app:trubin}

\begin{proof}[Proof of Proposition~\ref{prop:sparse-counterexample}]
We use the coordinates $(z_1,z_2,z_3,z_{12},z_{123})$, so that $z^S=(1,1,0,1,0)$ and $z^T=(0,0,1,0,0)$. We first show that $z^S$ and $z^T$ are adjacent in $\MP(H)$. For $c=(-1,-1,1,3,-2)$, direct enumeration gives the following table, where $z^R$ denotes the multilinear vertex of $R\subseteq V$.
\[
\begin{array}{c|cccccccc}
R & \emptyset & \{1\} & \{2\} & \{3\} & \{1,2\} & \{1,3\} & \{2,3\} & \{1,2,3\}\\ \hline
c^\top z^R & 0 & -1 & -1 & 1 & 1 & 0 & 0 & 0
\end{array}
\]
Hence $z^S$ and $z^T$ are the only vertices of $\MP(H)$ that maximize $c^\top z$, so that $\MP(H)\cap\{z:c^\top z=1\}=\conv\{z^S,z^T\}$, and $z^S$ and $z^T$ are adjacent in $\MP(H)$.

We now show that $z^S$ and $z^T$ are not adjacent in $\SR(H)$. We have $S\setminus T=\{1,2\}$, $T\setminus S=\{3\}$ and $S\cap T=\emptyset$, so $\Theta(H,S,T)$ has node set $\{1,2,3\}$. The edge $\{1,2\}$ has two nodes in $S\setminus T$ and none in $T\setminus S$, so it satisfies condition (a) of Definition~\ref{def:contracted-hypergraph} and contributes the edge $\{1,2\}$ to $\Theta(H,S,T)$. The edge $\{1,2,3\}$ has two nodes in $S\setminus T$ and one in $T\setminus S$, so it satisfies none of conditions (a)--(c) and contributes no edge. Therefore $E_\Theta=\bigl\{\{1,2\}\bigr\}$, and $\Theta(H,S,T)$ has the two connected components $\{1,2\}$ and $\{3\}$. By Theorem~\ref{thm:adjacency-sr}, $z^S$ and $z^T$ are not adjacent in $\SR(H)$.
\end{proof}

The hypergraph in Proposition~\ref{prop:sparse-counterexample} is as small as possible.
\begin{remark}\label{rem:minimality}
If $|E|\le1$, then $\SR(H)=\MP(H)$. Indeed, with no edges both polytopes equal $[0,1]^V$. With a single edge $e$, the inequalities of $\SR(H)$ for $e$ describe the convex hull of the binary solutions of $z_e=\prod_{v\in e}z_v$ \cite{Cra93}, and $\SR(H)$ is the product of this polytope with $[0,1]^{V\setminus e}$, as is $\MP(H)$. A hypergraph with at most two nodes has at most one edge. Hence Proposition~\ref{prop:sparse-counterexample} uses the minimum possible number of nodes and edges.
\end{remark}

We now describe the face responsible for the failure of the Trubin property in Proposition~\ref{prop:sparse-counterexample}.
\begin{remark}\label{rem:counterexample-face}
Only the constraints $z_{123}\ge 0$, $z_{12}\le z_1$ and $z_{12}\le z_2$ are tight at the midpoint $m = (\frac{1}{2},\frac{1}{2},\frac{1}{2},\frac{1}{2},0)$, so the minimal face of $\SR(H)$ containing $z^S$ and $z^T$ is
\[
    F^m = \{z\in \SR(H): z_{123} = 0,\ z_{12} = z_1 = z_2\}.
\]
In the parameters $t := z_1 = z_2 = z_{12}$ and $s := z_3$, the remaining inequalities of $\SR(H)$ read $0\le t\le 1$, $0\le s\le 1$ and $2t+s\le 2$, so that $F^m$ is the quadrilateral with vertices
\[
    z^{\emptyset} = (0,0,0,0,0),\quad z^{S} = (1,1,0,1,0),\quad \bar{z} = \Bigl(\tfrac{1}{2},\tfrac{1}{2},1,\tfrac{1}{2},0\Bigr),\quad z^{T} = (0,0,1,0,0),
\]
in cyclic order, shown in Fig.~\ref{fig:counterexample}(b). The only fractional vertex of this face is $\bar z$, and it is adjacent to both $z^S$ and $z^T$. Hence the segment $[z^S,z^T]$ is a diagonal of the quadrilateral. Intersecting $F^m$ with $\MP(H)$ gives instead the triangle $\conv\{z^{\emptyset},z^S,z^T\}$, in which $z^S$ and $z^T$ are adjacent. Thus the Trubin property fails because a fractional vertex of $\SR(H)$ lies in the minimal face. For graphs, Theorem~\ref{thm:cube} excludes such a fractional vertex.
\end{remark}

The following lemma shows that every failure of the Trubin property produces a fractional vertex, as in Remark~\ref{rem:counterexample-face}.
\begin{lemma}\label{lem:fractional-failure}
Let $H$ be a hypergraph, and let $S,T\subseteq V$ be such that $z^S$ and $z^T$ are adjacent in $\MP(H)$ but not in $\SR(H)$. Then the minimal face of $\SR(H)$ containing $z^S$ and $z^T$ has a fractional vertex.
\end{lemma}

\begin{proof}
Let $z^S$ and $z^T$ be adjacent in $\MP(H)$ but not in $\SR(H)$, and let $F$ be the minimal face of $\SR(H)$ containing them. If all vertices of $F$ were integral, then $F$ would be contained in $\MP(H)$. Since the inequalities of $\SR(H)$ are valid for $\MP(H)$, the polytope $F$ would then be a face of $\MP(H)$. The segment between $z^S$ and $z^T$, which is a face of $\MP(H)$ contained in $F$, would be a face of $F$, and hence of $\SR(H)$, a contradiction. Hence $F$ has a fractional vertex.
\end{proof}

\subsection{The minimal face as a cube}
\begin{proof}[Proof of Theorem~\ref{thm:cube}]
Let $m$, $U$, $\sigma$ and $\mathcal W^m$ be as in Appendix~\ref{sec:app-adjacency}, and let $F^m$ be the minimal face of $\SR(H)$ containing $m$, which is also the minimal face containing $z^S$ and $z^T$. By the preliminaries at the beginning of Appendix~\ref{app:adjacency}, $F^m=\SR(H)\cap(m+\mathcal W^m)$, and $m$ lies in the relative interior of $F^m$.
For $t\in\mathbb R^k$, define $w(t)\in\mathbb R^{V\cup E}$ by $w(t)_v=\sigma(v)t_l$ for $v\in C_l$, $w(t)_v=0$ for $v\notin U$, $w(t)_f=0$ if $m_f\in\{0,1\}$, and $w(t)_f=w(t)_v$ for any $v\in f\cap U$ if $m_f=\frac12$. The proof of Lemma~\ref{lem:projection} shows that $w(t)\in\mathcal W^m$, and the injectivity of $\pi$ shown in the proof of Proposition~\ref{prop:minimal-face-dimension} gives $\mathcal W^m=\{w(t):t\in\mathbb R^k\}$. The map $t\mapsto m+w(t)$ is affine and injective.

We show that $z:=m+w(t)$ lies in $\SR(H)$ if and only if $t\in[-\frac12,\frac12]^k$. For $v\in C_l$ we have $z_v=\frac12+\sigma(v)t_l$, so the node bounds of the nodes in $U$ hold if and only if $|t_l|\le\frac12$ for every $l$. Assume now $t\in[-\frac12,\frac12]^k$. Every node $v\notin U$ has $z_v=m_v\in\{0,1\}$. Let $f\in E$.
If $m_f=1$, then $f\subseteq S\cap T$, and $z_v=z_f=1$ for every $v\in f$, so all inequalities of $f$ hold.
If $m_f=\frac12$, then $f$ is contained in exactly one of $S$ and $T$. Then $f\setminus U\subseteq S\cap T$, and $f\cap U$ is a nonempty subset of $S\setminus T$ or of $T\setminus S$. If $|f\cap U|\ge2$, then $f$ satisfies condition (a) or (b) of Definition~\ref{def:contracted-hypergraph}, so $f\cap U\in E_\Theta$. In both cases $f\cap U\subseteq C_l$ for some $l$, and $\sigma$ is constant on $f\cap U$. Hence all nodes $v\in f\cap U$ have the same coordinate $z_v=\alpha:=\frac12+\sigma(v)t_l\in[0,1]$, and $z_f=\alpha$. Let $\gamma:=|f\cap U|\ge1$. Since $z_v=1$ for $v\in f\setminus U$, we have $0\le z_f\le z_v$ for every $v\in f$, and $\sum_{v\in f}z_v-|f|+1=1-\gamma(1-\alpha)\le\alpha=z_f$.
If $m_f=0$, then $z_f=0\le z_v$ for every $v\in f$, and it remains to show that $\sum_{v\in f}(1-z_v)\ge1$. If $f$ has a node $x\notin S\cup T$, then $z_x=0$. Otherwise $f\subseteq S\cup T$ is contained in neither $S$ nor $T$, so $f$ meets both $S\setminus T$ and $T\setminus S$. By $(\star)$, there are $u\in f\cap(S\setminus T)$ and $v\in f\cap(T\setminus S)$ in the same component $C_l$, and $(1-z_u)+(1-z_v)=(\frac12-t_l)+(\frac12+t_l)=1$. All other terms of the sum are nonnegative.

Therefore $F^m=\{m+w(t):t\in[-\frac12,\frac12]^k\}$ is the image of a cube under an injective affine map. Its vertices are the points $m+w(t)$ with $t\in\{-\frac12,\frac12\}^k$. For such $t$, every node coordinate is $0$ or $1$, and every edge coordinate is $m_f\in\{0,1\}$ or $\frac12+\sigma(v)t_l\in\{0,1\}$, so these vertices are integral. The node coordinates on $C_l$ agree with those of $z^S$ if $t_l=\frac12$ and with those of $z^T$ if $t_l=-\frac12$. Hence they agree with the node coordinates of $z^{S\triangle C_I}$, where $I=\{l:t_l=-\frac12\}$. Since an integral point of $\SR(H)$ is the multilinear point determined by its node coordinates, the vertex equals $z^{S\triangle C_I}$.

Finally, all vertices of $F^m$ lie in $\MP(H)$, so $F^m\subseteq\MP(H)$. Let $c^\top z\le\delta$ be an inequality valid for $\SR(H)$ with $F^m=\SR(H)\cap\{z:c^\top z=\delta\}$. It is also valid for $\MP(H)\subseteq\SR(H)$, so $F^m=\MP(H)\cap\{z:c^\top z=\delta\}$ is a face of $\MP(H)$. Let $F'$ be a face of $\MP(H)$ containing $z^S$ and $z^T$. Then $F'\cap F^m$ is a face of $F^m$ containing $m$. Since $m$ lies in the relative interior of $F^m$, we have $F'\cap F^m=F^m$. Hence $F^m$ is the minimal face of $\MP(H)$ containing $z^S$ and $z^T$. Its dimension is $k$, so $z^S$ and $z^T$ are adjacent in $\MP(H)$ if and only if $k=1$.
\end{proof}

\subsection{Local and structural sufficient conditions}
Theorem~\ref{thm:cube} gives the following sufficient condition for a single pair $S,T$.
\begin{corollary}\label{cor:local-trubin}
    Let $H$ be a hypergraph and let $S,T\subseteq V$ be distinct. Set $U:= S\triangle T$ and let $C_1,\ldots,C_k$ be the connected components of $\Theta(H,S,T)$. Suppose that for every edge $e\in E$ with $e\subseteq S\cup T$, there exists $l\in [k]$ such that $e\cap U \subseteq C_l$. Then $z^S$ and $z^T$ are adjacent in $\MP(H)$ if and only if $\Theta(H,S,T)$ is connected.
\end{corollary}

\begin{proof}
Let $e\subseteq S\cup T$ be an edge that meets both $S\setminus T$ and $T\setminus S$. By hypothesis, $e\cap U$ lies in a single connected component of $\Theta(H,S,T)$, so any $u\in e\cap(S\setminus T)$ and $v\in e\cap(T\setminus S)$ lie in the same component. Hence $S,T$ satisfy $(\star)$, and Theorem~\ref{thm:cube} gives the statement.
\end{proof}

\begin{proof}[Proof of Theorem~\ref{thm:matching-trubin}]
Fix distinct $S,T\subseteq V$, let $U:=S\triangle T$, and let $C_1,\ldots,C_k$ be the connected components of $\Theta(H,S,T)$. By Corollary~\ref{cor:local-trubin}, it suffices to show that $e\cap U$ lies in a single component for every edge $e\in E$ with $e\subseteq S\cup T$. If $|e\cap U|\le1$, there is nothing to prove. If $e\cap U=\{u,v\}$, then $e\setminus U\subseteq S\cap T$, so $e$ satisfies one of conditions (a)--(c) of Definition~\ref{def:contracted-hypergraph}. Hence $\{u,v\}\in E_\Theta$. Now suppose that $|e\cap U|\ge3$, and let $u,v\in e\cap U$ be distinct. If $\{u,v\}\in E$, then $\{u,v\}$ satisfies one of conditions (a)--(c), so $\{u,v\}\in E_\Theta$. Otherwise, $\{u,v\}\in M_e$. Let $r\in(e\cap U)\setminus\{u,v\}$. Since $M_e$ is a matching, neither $\{u,r\}$ nor $\{v,r\}$ belongs to $M_e$. Hence both are edges of $H$, and therefore edges of $\Theta(H,S,T)$. In all cases $u$ and $v$ lie in the same component. Corollary~\ref{cor:local-trubin} now gives the first statement, and the Trubin property follows by Theorem~\ref{thm:adjacency-sr}.
\end{proof}

\subsection{Bad pairs and acyclic traces}
Throughout this subsection, let $S,T\subseteq V$ be distinct, let $U:=S\triangle T$, and let $X:=S\cap T$. Recall that the trace of an edge $e\subseteq S\cup T$ is $e\cap U$, and that $e\setminus U\subseteq X$. Every edge $e\subseteq S\cup T$ whose trace has exactly two nodes satisfies one of conditions (a)--(c) of Definition~\ref{def:contracted-hypergraph}. Every edge $e\subseteq S\cup T$ whose trace has at least two nodes and is contained in $S\setminus T$ or in $T\setminus S$ satisfies condition (a) or (b). In both cases $e$ is a generating edge, so its trace lies in a single connected component of $\Theta(H,S,T)$.

\begin{proof}[Proof of Theorem~\ref{thm:bad-pair}]
Let $S':=X\cup K$ and $T':=X\cup K'$, so that $S'\setminus T'=K$ and $T'\setminus S'=K'$. We first show that $\Theta(H,S',T')$ has exactly the two connected components $K$ and $K'$. Let $f$ be a generating edge for the pair $S,T$ whose trace $f\cap U$ lies in $K$ and has at least two nodes. Since $f\setminus U\subseteq X$, we have $f\subseteq X\cup K$, so $f$ satisfies condition (a) for the pair $S',T'$, with the same trace. Since $K$ is connected in $\Theta(H,S,T)$, it is therefore connected in $\Theta(H,S',T')$. Similarly, condition (b) shows that $K'$ is connected in $\Theta(H,S',T')$. Suppose that some generating edge $f$ of $\Theta(H,S',T')$ meets both $K$ and $K'$. Then $f$ satisfies condition (c) for $S',T'$, so $f\subseteq X\cup K\cup K'\subseteq S\cup T$, and the trace of $f$ with respect to $S,T$ consists of one node of $K$ and one node of $K'$. Hence this trace lies in a single component of $\Theta(H,S,T)$, a contradiction. Thus $\Theta(H,S',T')$ has the two connected components $K$ and $K'$, and $z^{S'}$ and $z^{T'}$ are not adjacent in $\SR(H)$ by Theorem~\ref{thm:adjacency-sr}.

We now show that $z^{S'}$ and $z^{T'}$ are adjacent in $\MP(H)$. Let $m':=\frac12(z^{S'}+z^{T'})$, and suppose that $m'=\sum_R\lambda_Rz^R$, where the sum is over $R\subseteq V$, $\lambda\ge0$, and $\sum_R\lambda_R=1$. Let $R$ satisfy $\lambda_R>0$. Since $m'_v=1$ for $v\in X$ and $m'_v=0$ for $v\notin X\cup K\cup K'$, we have $X\subseteq R\subseteq X\cup K\cup K'$. Let $f$ be a generating edge for the pair $S,T$ whose trace lies in $K$ and has at least two nodes, and let $v\in f\cap K$. Then $f\subseteq S'$ and $f\not\subseteq T'$, so $m'_f=\frac12=m'_v$. Since $z^R_f\le z^R_v$ for every $R$, the equation $\sum_R\lambda_R z^R_f=\sum_R\lambda_R z^R_v$ forces $z^R_f=z^R_v$ whenever $\lambda_R>0$. Hence $v\in R$ implies $f\subseteq R$. Since these traces connect $K$, we obtain $R\cap K\in\{\emptyset,K\}$, and similarly $R\cap K'\in\{\emptyset,K'\}$. Moreover, the edge $e$ satisfies $e\subseteq X\cup K\cup K'$ and meets both $K$ and $K'$. Hence $e$ is contained in neither $S'$ nor $T'$, so $m'_e=0$, and $R$ does not contain $K\cup K'$. Therefore $R\in\{X,X\cup K,X\cup K'\}$. Comparing the coordinates of a node of $K$ and of a node of $K'$ gives $\lambda_{X\cup K}=\lambda_{X\cup K'}=\frac12$, so $\lambda_X=0$. Thus $m'=\frac12z^{S'}+\frac12z^{T'}$ is the only way to write $m'$ as a convex combination of vertices of $\MP(H)$. If $z^{S'}$ and $z^{T'}$ were not adjacent in $\MP(H)$, the minimal face of $\MP(H)$ containing them would have dimension at least two. Since $m'$ lies in the relative interior of this face, it would be a convex combination with positive weight on some other vertex of the face, a contradiction. Hence $z^{S'}$ and $z^{T'}$ are adjacent in $\MP(H)$, and $\SR(H)$ does not have the Trubin property with respect to $\MP(H)$.
\end{proof}

A connected component $K$ of $\Theta(H,S,T)$ is \emph{blocked} by an edge $e\subseteq S\cup T$ with trace $t$ if $t$ meets both $K$ and $U\setminus K$, and either $t\cap K\subseteq S\setminus T$ and $t\setminus K\subseteq T\setminus S$, or $t\cap K\subseteq T\setminus S$ and $t\setminus K\subseteq S\setminus T$. Whether $K$ is blocked by $e$ depends only on the trace of $e$.

\begin{lemma}\label{lem:free-component}
Suppose that $\Theta(H,S,T)$ has at least two connected components, and that some component $K$ is blocked by no edge. Then $z^S$ and $z^T$ are not adjacent in $\MP(H)$.
\end{lemma}
\begin{proof}
Let $R:=S\triangle K$ and $R':=T\triangle K$. We claim that $z^S+z^T=z^R+z^{R'}$. For $v\notin U$, all four vertices have the same coordinate $z_v$. For $v\in U\setminus K$, we have $z^R_v=z^S_v$ and $z^{R'}_v=z^T_v$, and for $v\in K$, we have $z^R_v=z^T_v$ and $z^{R'}_v=z^S_v$. Let $e\in E$. Since $R\cup R'=S\cup T$, all four vertices have $z_e=0$ if $e\not\subseteq S\cup T$. Otherwise, let $t$ be the trace of $e$. If $t\cap K=\emptyset$, then $z^R_e=z^S_e$ and $z^{R'}_e=z^T_e$. If $t\subseteq K$, then $z^R_e=z^T_e$ and $z^{R'}_e=z^S_e$. Suppose that $t$ meets both $K$ and $U\setminus K$. Since $t$ does not lie in a single component, $t$ is contained neither in $S\setminus T$ nor in $T\setminus S$. Hence $e$ is contained neither in $S$ nor in $T$. Moreover, $e\subseteq R$ if and only if $t\cap K\subseteq T\setminus S$ and $t\setminus K\subseteq S\setminus T$, and $e\subseteq R'$ if and only if $t\cap K\subseteq S\setminus T$ and $t\setminus K\subseteq T\setminus S$. Since $K$ is not blocked by $e$, neither holds, so all four vertices have $z_e=0$. This proves the claim. Since $\emptyset\neq K\neq U$, the set $R$ differs from both $S$ and $T$, so $\{z^R,z^{R'}\}\neq\{z^S,z^T\}$. By Lemma~\ref{lem:midpoint-nonadjacency}, $z^S$ and $z^T$ are not adjacent in $\MP(H)$.
\end{proof}

If every trace lies in a single component, then no component is blocked, and Lemma~\ref{lem:free-component} gives another proof of the nontrivial direction of Corollary~\ref{cor:local-trubin}.

\begin{lemma}\label{lem:blocked-count}
Suppose that the pair $S,T$ is not bad. Then every trace blocks at most one component. If moreover $\mathcal G(S,T)$ is a forest and $\Theta(H,S,T)$ has $k\ge2$ components, then at most $(k-1)/2$ components are blocked.
\end{lemma}
\begin{proof}
A trace that meets at most one component blocks no component. Since the pair $S,T$ is not bad, no trace meets exactly two components. So consider a trace $t$ that meets at least three components. Suppose that $t$ blocks two distinct components $K$ and $K'$. By symmetry, we may assume $t\cap K\subseteq S\setminus T$ and $t\setminus K\subseteq T\setminus S$. Then $t\cap K'\subseteq T\setminus S$, so $t\setminus K'\subseteq S\setminus T$. Let $M$ be a third component met by $t$. Then $t\cap M$ is contained in both $t\setminus K\subseteq T\setminus S$ and $t\setminus K'\subseteq S\setminus T$, so $t\cap M=\emptyset$, a contradiction. Hence $t$ blocks at most one component.

Now suppose that $\mathcal G(S,T)$ is a forest, and let $r$ be the number of distinct traces that meet at least two components. Each of them meets at least three components, so $\mathcal G(S,T)$ has at least $3r$ edges and $k+r$ nodes. Since a forest has fewer edges than nodes, we get $3r\le k+r-1$, that is, $r\le(k-1)/2$. Since blocking depends only on the trace and each trace blocks at most one component, at most $r$ components are blocked.
\end{proof}

\begin{proof}[Proof of Theorem~\ref{thm:acyclic-traces}]
If $\Theta(H,S,T)$ is connected, then $z^S$ and $z^T$ are adjacent in $\SR(H)$ by Theorem~\ref{thm:adjacency-sr}, and hence in $\MP(H)$ by Lemma~\ref{lem:adjacency-subpolytope}. Conversely, suppose that $\Theta(H,S,T)$ has $k\ge2$ components. By Lemma~\ref{lem:blocked-count}, at most $(k-1)/2<k$ components are blocked, so some component is blocked by no edge. By Lemma~\ref{lem:free-component}, $z^S$ and $z^T$ are not adjacent in $\MP(H)$.
\end{proof}

\subsection{Hypergraphs of rank at most three}\label{sec:rank-three}
In this subsection we prove Theorem~\ref{thm:rank-three}.

We keep the notation $S$, $T$, $U$ and $X$ of the previous subsection. We use the following three facts.
\begin{itemize}
\item[(F1)] Let $f$ be a generating edge. Then its trace $f\cap U$ is an edge of $\Theta(H,S,T)$, it has at least two nodes, and it lies in a single connected component of $\Theta(H,S,T)$. Indeed, each of conditions (a)--(c) of Definition~\ref{def:contracted-hypergraph} gives $|f\cap U|\ge2$.
\item[(F2)] Every edge $e\subseteq S\cup T$ whose trace has exactly two nodes is a generating edge. Indeed, if both nodes of the trace lie in $S\setminus T$, then $e$ satisfies condition (a). If both lie in $T\setminus S$, then $e$ satisfies condition (b). Otherwise $e$ satisfies condition (c).
\item[(F3)] Suppose that $H$ has rank at most three and that the pair $S,T$ is not bad. Let $t$ be a trace that meets at least two connected components of $\Theta(H,S,T)$. Then $t$ meets exactly three connected components and contains exactly one node of each of them. Moreover, every edge $e\subseteq S\cup T$ with trace $t$ satisfies $e=t$. Indeed, since the pair is not bad, $t$ meets at least three components. Since $t\subseteq e$ and $|e|\le3$, the set $t$ has exactly three nodes, one in each of these components, and $e=t$.
\end{itemize}
The first lemma puts one connected component $L$ of $\Theta(H,S,T)$ alone on one side of the symmetric difference, and all other nodes of $U$ on the other side. It holds for hypergraphs of arbitrary rank. Since all nodes of $U\setminus L$ are then on the same side, the connected components of $\Theta(H,S,T)$ other than $L$ can only merge. The set $L$ remains a connected component.
\begin{lemma}\label{lem:one-sided}
Let $H$ be a hypergraph, and let $L$ be a connected component of $\Theta(H,S,T)$. Let $S':=(S\cup T)\setminus L$ and $T':=X\cup L$. Then the following statements hold.
\begin{itemize}
\item[(i)] Every edge $e\subseteq S'$ whose trace $e\cap U$ has at least two nodes is a generating edge for the pair $S',T'$, and its trace with respect to $S',T'$ is $e\cap U$.
\item[(ii)] Every connected component of $\Theta(H,S,T)$ lies in a single connected component of $\Theta(H,S',T')$.
\item[(iii)] The set $L$ is a connected component of $\Theta(H,S',T')$.
\end{itemize}
\end{lemma}
\begin{proof}
We have $S'\cap T'=X$, $S'\setminus T'=U\setminus L$ and $T'\setminus S'=L$. Hence $S'\triangle T'=U$ and $S'\cup T'=S\cup T$. In particular, every edge contained in $S\cup T$ has the same trace with respect to $S,T$ and with respect to $S',T'$.

(i) Let $e\subseteq S'$ satisfy $|e\cap U|\ge2$. Then $e\cap(T'\setminus S')=e\cap L=\emptyset$ and $|e\cap(S'\setminus T')|=|e\cap U|\ge2$. Hence $e$ satisfies condition (a) for $S',T'$, and its trace is $e\cap U$.

(ii) Let $C$ be a connected component of $\Theta(H,S,T)$. We show that every edge $g$ of $\Theta(H,S,T)$ contained in $C$ is also an edge of $\Theta(H,S',T')$. Let $f$ be a generating edge for the pair $S,T$ with $f\cap U=g$. By (F1), $g$ has at least two nodes. Since $f\setminus U\subseteq X$, we have $f\subseteq X\cup C$. If $C\neq L$, then $C\subseteq U\setminus L$, so $f\subseteq S'$, and $g$ is an edge of $\Theta(H,S',T')$ by (i). If $C=L$, then $f\subseteq T'$, $f\cap(S'\setminus T')=\emptyset$ and $|f\cap(T'\setminus S')|=|g|\ge2$. Hence $f$ satisfies condition (b) for $S',T'$, and $g$ is again an edge of $\Theta(H,S',T')$. Any two nodes of $C$ are joined by a path in $\Theta(H,S,T)$ whose edges are contained in $C$. This path is also a path in $\Theta(H,S',T')$, which proves (ii).

(iii) Suppose that some edge of $\Theta(H,S',T')$ meets both $L$ and $U\setminus L$, and let $f$ be a generating edge for the pair $S',T'$ whose trace is this edge. Then $f$ meets both $S'\setminus T'=U\setminus L$ and $T'\setminus S'=L$, so $f$ satisfies condition (c) for $S',T'$. Hence the trace of $f$ consists of a node $u\in U\setminus L$ and a node $v\in L$. Since $f\subseteq S'\cup T'=S\cup T$, the edge $f$ is a generating edge for the pair $S,T$ by (F2). By (F1), $u$ and $v$ lie in the same connected component of $\Theta(H,S,T)$, which contradicts $u\notin L$ and $v\in L$. Therefore every edge of $\Theta(H,S',T')$ is contained in $L$ or in $U\setminus L$. Together with (ii), this shows that $L$ is a connected component of $\Theta(H,S',T')$.
\end{proof}

The second lemma treats two traces that meet the same three connected components.
\begin{lemma}\label{lem:same-components}
Let $H$ be a hypergraph of rank at most three, and suppose that the pair $S,T$ is not bad. Let $e\neq e'$ be edges contained in $S\cup T$ whose traces meet exactly the same three connected components $K_1,K_2,K_3$ of $\Theta(H,S,T)$. Then $e\setminus e'\neq\emptyset$. Moreover, for every $p\in e\setminus e'$, the pair $S\cup\{p\},T\cup\{p\}$ is bad.
\end{lemma}
\begin{proof}
By (F3), we have $e=e\cap U$ and $e'=e'\cap U$, and each of $e$ and $e'$ consists of one node of each of $K_1$, $K_2$ and $K_3$. Since $|e|=|e'|=3$ and $e\neq e'$, the set $e\setminus e'$ is nonempty. Let $p\in e\setminus e'$, and assume without loss of generality that $p\in K_1$. Let $q$ be the node of $e'$ in $K_1$. Since $p\notin e'$, we have $q\neq p$. Let $S':=S\cup\{p\}$ and $T':=T\cup\{p\}$. Then $S'\cup T'=S\cup T$, $S'\cap T'=X\cup\{p\}$, $S'\setminus T'=(S\setminus T)\setminus\{p\}$, $T'\setminus S'=(T\setminus S)\setminus\{p\}$ and $U':=S'\triangle T'=U\setminus\{p\}$. We prove three claims.

\textit{Claim A.} Every generating edge $f$ of $\Theta(H,S,T)$ with $p\notin f$ is a generating edge for the pair $S',T'$ with the same trace. Indeed, $f\subseteq S'\cup T'$, $f\cap(S'\setminus T')=f\cap(S\setminus T)$, $f\cap(T'\setminus S')=f\cap(T\setminus S)$ and $f\cap U'=f\cap U$.

\textit{Claim B.} The edge $e$ is a generating edge for the pair $S',T'$ with trace $e\setminus\{p\}$. Indeed, $e\subseteq S'\cup T'$, and its trace with respect to $S',T'$ is $e\cap U'=e\setminus\{p\}$, which has exactly two nodes. The claim follows from (F2), applied to the pair $S',T'$ in place of $S,T$.

\textit{Claim C.} No edge of $\Theta(H,S',T')$ meets both $K_1\setminus\{p\}$ and $U\setminus K_1$. Suppose otherwise, and let $f$ be a generating edge for the pair $S',T'$ whose trace $f\cap U'$ contains a node $u\in K_1\setminus\{p\}$ and a node $v\in U\setminus K_1$. Since $f\subseteq S'\cup T'=S\cup T$, the trace $f\cap U$ of $f$ with respect to $S,T$ is defined, and it meets $K_1$ and the connected component of $\Theta(H,S,T)$ containing $v$. By (F3), we have $f=f\cap U$, and $f$ consists of one node of each of three distinct connected components of $\Theta(H,S,T)$. Its node in $K_1$ is $u$, so $p\notin f$. Hence $f\cap(S'\setminus T')=f\cap(S\setminus T)$ and $f\cap(T'\setminus S')=f\cap(T\setminus S)$. Since $f$ satisfies one of conditions (a)--(c) of Definition~\ref{def:contracted-hypergraph} for $S',T'$, it satisfies the same condition for $S,T$, and $f$ is a generating edge for the pair $S,T$. By (F1), its trace lies in a single connected component of $\Theta(H,S,T)$, a contradiction.

We now conclude. For $j\in\{2,3\}$, every generating edge $f$ of $\Theta(H,S,T)$ whose trace lies in $K_j$ satisfies $f\subseteq X\cup K_j$, and in particular $p\notin f$. By Claim A, every edge of $\Theta(H,S,T)$ contained in $K_2$ or in $K_3$ is an edge of $\Theta(H,S',T')$. Hence the nodes of $K_2$ lie in a single connected component of $\Theta(H,S',T')$, and so do the nodes of $K_3$. By Claim B, the node of $e$ in $K_2$ and the node of $e$ in $K_3$ form an edge of $\Theta(H,S',T')$. Hence $K_2\cup K_3$ is contained in a single connected component $M$ of $\Theta(H,S',T')$. By Claim C, every edge of $\Theta(H,S',T')$ that contains a node of $K_1\setminus\{p\}$ is contained in $K_1\setminus\{p\}$. Hence the connected component of $\Theta(H,S',T')$ containing $q$ is contained in $K_1\setminus\{p\}$, and it differs from $M$. Finally, $e'\subseteq S'\cup T'$ and $p\notin e'$, so the trace of $e'$ with respect to $S',T'$ is $e'\cap U'=e'$. It consists of $q$, a node of $K_2$ and a node of $K_3$. Hence it meets exactly two connected components of $\Theta(H,S',T')$, namely the one containing $q$ and $M$. Therefore the pair $S',T'$ is bad.
\end{proof}

We now show that the second hypothesis of Theorem~\ref{thm:acyclic-traces} holds for all pairs when $H$ has rank at most three and no bad pair. We write a cycle of length $2m$ in the bipartite graph $\mathcal G(S,T)$ as a sequence $K_1,t_1,K_2,t_2,\ldots,K_m,t_m$, where $m\ge2$, the sets $K_1,\ldots,K_m$ are distinct connected components of $\Theta(H,S,T)$, the sets $t_1,\ldots,t_m$ are distinct traces that meet at least two connected components, and $t_i$ meets $K_i$ and $K_{i+1}$ for every $i\in[m]$, where $K_{m+1}:=K_1$. The proof takes a shortest cycle of $\mathcal G(S,T)$ and constructs a bad pair. If the cycle consists of two traces that meet the same three connected components, we use Lemma~\ref{lem:same-components}. Otherwise, we use Lemma~\ref{lem:one-sided}.
\begin{proposition}\label{prop:rank-three-forest}
Let $H$ be a hypergraph of rank at most three with no bad pair, and let $S,T\subseteq V$ be distinct. Then $\mathcal G(S,T)$ is a forest.
\end{proposition}
\begin{proof}
Suppose that $\mathcal G(S,T)$ contains a cycle, and let $K_1,t_1,\ldots,K_m,t_m$ be a cycle of minimum length $2m$. Since the pair $S,T$ is not bad, (F3) shows that each trace $t_i$ meets exactly three connected components and is itself an edge of $H$. For $i\in[m]$, let $L_i$ be the connected component met by $t_i$ other than $K_i$ and $K_{i+1}$.

Suppose first that $m=2$ and $L_1=L_2$. Then $t_1$ and $t_2$ are distinct edges contained in $S\cup T$ whose traces meet exactly the same three connected components $K_1$, $K_2$ and $L_1$. By Lemma~\ref{lem:same-components}, the pair $S\cup\{p\},T\cup\{p\}$ is bad for every $p\in t_1\setminus t_2$, a contradiction.

Suppose now that $m\ge3$, or that $m=2$ and $L_1\neq L_2$. We claim that none of $t_1,\ldots,t_{m-1}$ meets $L_m$. Since $t_m$ meets three distinct connected components, we have $L_m\notin\{K_m,K_1\}$. If $m=2$, then $t_1$ meets exactly $K_1$, $K_2$ and $L_1$, and $L_2$ is none of them. Let $m\ge3$, and suppose that $t_i$ meets $L_m$ for some $i\in[m-1]$. If $L_m=K_j$ for some $j$ with $2\le j\le m-1$, then $K_1,t_1,\ldots,K_{j-1},t_{j-1},K_j,t_m$ is a cycle of $\mathcal G(S,T)$ of length $2j<2m$. Otherwise $L_m$ is not one of $K_1,\ldots,K_m$. Then
\[
L_m,t_i,K_{i+1},t_{i+1},\ldots,K_m,t_m
\qquad\text{and}\qquad
L_m,t_m,K_1,t_1,\ldots,K_i,t_i
\]
are cycles of $\mathcal G(S,T)$ of lengths $2(m-i+1)$ and $2(i+1)$. Their lengths sum to $2m+4$, so one of them has length at most $m+2<2m$. In both cases we contradict the minimality of $m$, which proves the claim.

Let $S':=(S\cup T)\setminus L_m$ and $T':=X\cup L_m$. Let $i\in[m-1]$. Since $t_i$ is an edge of $H$ contained in $U$ and disjoint from $L_m$, we have $t_i\subseteq S'$. Its trace $t_i$ has three nodes, so $t_i$ is an edge of $\Theta(H,S',T')$ by Lemma~\ref{lem:one-sided}(i). Each $K_j$ lies in a single connected component of $\Theta(H,S',T')$ by Lemma~\ref{lem:one-sided}(ii), and $t_i$ meets both $K_i$ and $K_{i+1}$. Hence $K_1\cup\cdots\cup K_m$ is contained in a single connected component $M$ of $\Theta(H,S',T')$. By Lemma~\ref{lem:one-sided}(iii), $L_m$ is a connected component of $\Theta(H,S',T')$, and $L_m\neq M$, because $K_1\subseteq M$ is disjoint from $L_m$. The edge $t_m$ is contained in $S\cup T=S'\cup T'$, and its trace with respect to $S',T'$ is $t_m\cap U=t_m$. It meets $K_m\cup K_1\subseteq M$ and $L_m$, and no other connected component of $\Theta(H,S',T')$. Hence the pair $S',T'$ is bad, a contradiction.
\end{proof}

\begin{proof}[Proof of Theorem~\ref{thm:rank-three}]
If $H$ has a bad pair, then $\SR(H)$ does not have the Trubin property with respect to $\MP(H)$ by Theorem~\ref{thm:bad-pair}. Conversely, suppose that $H$ has no bad pair, and let $S,T\subseteq V$ be distinct sets such that $z^S$ and $z^T$ are adjacent in $\MP(H)$. Since the pair $S,T$ is not bad, no trace meets exactly two connected components of $\Theta(H,S,T)$. By Proposition~\ref{prop:rank-three-forest}, $\mathcal G(S,T)$ is a forest. By Theorem~\ref{thm:acyclic-traces}, $\Theta(H,S,T)$ is connected. By Theorem~\ref{thm:adjacency-sr}, $z^S$ and $z^T$ are adjacent in $\SR(H)$. Since the vertices of $\MP(H)$ are exactly the multilinear vertices, the one-skeleton of $\MP(H)$ is contained in the one-skeleton of $\SR(H)$.
\end{proof}

The proof of Theorem~\ref{thm:rank-three} also bounds the number of blocked components. Suppose that $H$ has rank at most three and no bad pair, and that $\Theta(H,S,T)$ has $k\ge2$ connected components. By Proposition~\ref{prop:rank-three-forest} and Lemma~\ref{lem:blocked-count}, at most $(k-1)/2$ of them are blocked. By the proof of Lemma~\ref{lem:free-component}, switching any one of the at least $(k+1)/2$ remaining components in both $z^S$ and $z^T$ yields two further multilinear vertices with the same midpoint.

The proof is also constructive. Suppose that $H$ has rank at most three, and that $z^S$ and $z^T$ are adjacent in $\MP(H)$ while $\Theta(H,S,T)$ is disconnected. Then either some trace meets exactly two connected components, so that $S,T$ is itself a bad pair, or $\mathcal G(S,T)$ contains a cycle by Theorem~\ref{thm:acyclic-traces}. In the second case, the proof of Proposition~\ref{prop:rank-three-forest} applied to a shortest cycle yields the explicit bad pair $S\cup\{p\},T\cup\{p\}$ or $(S\cup T)\setminus L_m,X\cup L_m$. The rank assumption enters only through (F3), which is used in Lemma~\ref{lem:same-components} and in Proposition~\ref{prop:rank-three-forest}.

The sufficient conditions for the preservation of adjacency obtained so far can be ordered from the strongest to the weakest, except that the two weakest ones are incomparable:
\begin{align*}
\tbinom{V}{2}\subseteq E\ &\Longrightarrow\ \text{matching condition}\ \Longrightarrow\ \text{every trace lies in a single component}\\
&\Longrightarrow\ (\star)\ \text{and the hypotheses of Theorem~\ref{thm:acyclic-traces}}.
\end{align*}
If $\binom{V}{2}\subseteq E$, as in Theorem~\ref{thm:complete-two-edge-layer}, then $M_e=\emptyset$ for every $e\in E$, so the matching condition holds. The proof of Theorem~\ref{thm:matching-trubin} shows that the matching condition implies, for every pair $S,T$, the hypothesis of Corollary~\ref{cor:local-trubin}. This hypothesis in turn implies condition $(\star)$ of Theorem~\ref{thm:cube}, as shown in the proof of Corollary~\ref{cor:local-trubin}. Corollary~\ref{cor:local-trubin} is also the special case of Theorem~\ref{thm:acyclic-traces} in which every trace lies in a single component. Condition $(\star)$ and the hypotheses of Theorem~\ref{thm:acyclic-traces} are incomparable. For $E=\{\{1,4\},\{2,4\},\{3,4\},\{1,2,3\}\}$, $S=\{1\}$, and $T=\{2,3\}$, the contracted hypergraph has the components $\{1\}$, $\{2\}$, $\{3\}$, so $(\star)$ fails for the edge $\{1,2,3\}$. The hypotheses of Theorem~\ref{thm:acyclic-traces} hold for this pair. For $E=\{\{1,2\},\{1,2,3\}\}$, $S=\{1\}$, and $T=\{2,3\}$, condition $(\star)$ holds, but the trace $\{1,2,3\}$ meets exactly the two components $\{1,2\}$ and $\{3\}$. The matching condition and $\binom{V}{2}\subseteq E$ depend on $H$ alone, whereas the other conditions depend on the pair $S,T$. Hence Theorems~\ref{thm:cube} and~\ref{thm:acyclic-traces} give the two weakest sufficient conditions for the preservation of adjacency. Theorem~\ref{thm:matching-trubin} gives a condition on $H$ alone that guarantees the Trubin property. If $H$ is Berge-acyclic, then $\SR(H)=\MP(H)$ \cite{dPKha18SIOPT,BucCraRod19}, so $\SR(H)$ has the Trubin property with respect to $\MP(H)$. Berge-acyclicity and the matching condition are incomparable. A graph with a cycle satisfies the matching condition but is not Berge-acyclic, while the hypergraph with the single edge $\{1,2,3\}$ is Berge-acyclic but violates the matching condition.

\section{Diameter Bounds}\label{app:diameter}

\begin{proof}[Proof of Theorem~\ref{thm:diameter-mp}]
Let $C_1,\ldots,C_k$ be the connected components of $\Theta(H,S,T)$, let $S_0:=S$, and let $S_i:=S_{i-1}\triangle C_i$ for $i\in[k]$. Since $C_1,\ldots,C_k$ are pairwise disjoint and their union is $S\triangle T$, we have $S_k=S\triangle(S\triangle T)=T$. Fix $i\in[k]$. Then $S_{i-1}\triangle S_i=C_i$, and we show that the edges of $\Theta(H,S,T)$ contained in $C_i$ are also edges of $\Theta(H,S_{i-1},S_i)$. Let $g$ be such an edge, and let $f\in F_\Theta$ with $g=f\cap(S\triangle T)$. Then $f\cap C_j=\emptyset$ for $j\neq i$, and $f\cap(S\triangle T)=f\cap C_i$. Since intersection distributes over symmetric difference, we have $f\cap S_{i-1}=f\cap S$ and $f\cap S_i=(f\cap S)\triangle(f\cap C_i)=f\cap T$. Hence $f\cap(S_{i-1}\setminus S_i)=f\cap(S\setminus T)$ and $f\cap(S_i\setminus S_{i-1})=f\cap(T\setminus S)$. Moreover, $f\subseteq S\cup T$ gives $f\subseteq S_{i-1}\cup S_i$. Therefore $f$ is also a generating edge for the pair $S_{i-1},S_i$, and its trace is $f\cap C_i=g$.

Since $C_i$ is a connected component of $\Theta(H,S,T)$, the edges of $\Theta(H,S,T)$ contained in $C_i$ connect all nodes of $C_i$. Hence $\Theta(H,S_{i-1},S_i)$, whose node set is $C_i$, is connected. By Theorem~\ref{thm:adjacency-sr}, $z^{S_{i-1}}$ and $z^{S_i}$ are adjacent in $\SR(H)$, and by Lemma~\ref{lem:adjacency-subpolytope} they are also adjacent in $\MP(H)$. Hence $\dist_{\SR(H)}(z^S,z^T)\le k$ and $\dist_{\MP(H)}(z^S,z^T)\le k$. Since $C_1,\ldots,C_k$ are nonempty and pairwise disjoint subsets of $S\triangle T$, we have $k\le|S\triangle T|\le|V|$. Taking the maximum over all distinct $S,T\subseteq V$ yields $\diam(\MP(H))\le|V|$.
\end{proof}

Theorem~\ref{thm:complete-two-edge-layer} also characterizes the hypergraphs $H$ for which $\MP(H)$ has diameter one.
\begin{corollary}\label{cor:diameter-one-mp}
    Let $H$ be a hypergraph with $V\neq\emptyset$. Then $\diam(\MP(H)) = 1$ if and only if $\binom{V}{2} \subseteq E$.
\end{corollary}
\begin{proof}
    First suppose $\binom{V}{2} \subseteq E$. By Theorem~\ref{thm:complete-two-edge-layer}, every pair of distinct multilinear vertices $z^S$ and $z^T$ is adjacent in $\SR(H)$. Hence, by Lemma~\ref{lem:adjacency-subpolytope}, they are adjacent in $\MP(H)$ as well. As $V\neq\emptyset$, the polytope $\MP(H)$ has at least two vertices, and therefore $\diam(\MP(H)) = 1$.

    Conversely, suppose that $\{u,v\}\notin E$ for some pair of distinct nodes $u,v\in V$, and consider the four distinct multilinear vertices $z^{\emptyset}$, $z^{\{u,v\}}$, $z^{\{u\}}$, and $z^{\{v\}}$. Every edge of $H$ has cardinality at least two and is different from $\{u,v\}$, hence it is contained in none of the four sets $\emptyset$, $\{u,v\}$, $\{u\}$, $\{v\}$. Therefore all the edge coordinates of these four vertices vanish, and comparing the node coordinates gives
    \[
        z^{\emptyset} + z^{\{u,v\}} = z^{\{u\}} + z^{\{v\}}.
    \]
    Since $\{z^{\emptyset},z^{\{u,v\}}\}\neq\{z^{\{u\}},z^{\{v\}}\}$, Lemma~\ref{lem:midpoint-nonadjacency} shows that $z^{\emptyset}$ and $z^{\{u,v\}}$ are not adjacent in $\MP(H)$. Hence $\diam(\MP(H))\ge 2$.
\end{proof}

For graphs, Corollary~\ref{cor:diameter-one-mp} reads $\diam(\MP(G)) = 1$ if and only if $G$ is complete, which is due to Padberg \cite[Theorem~2]{Pad89}. Padberg observed that $\diam(\MP(G))\ge 2$ whenever $G$ is disconnected \cite[p.~157]{Pad89}. By Corollary~\ref{cor:diameter-one-mp}, the same holds for every connected graph $G$ that is not complete.

\subsection{Fractional vertices of graphs}
For a graph $G=(V,E)$ and $U\subseteq V$, we denote by $G-U$ the graph obtained from $G$ by deleting the nodes in $U$ and the edges incident to them.

\begin{lemma}\label{lem:graph-face}
Let $G=(V,E)$ be a graph, let $U\subseteq V$, and let $c\in\{0,1\}^U$. Then $F:=\{z\in\SR(G): z_v=c_v\ \forall v\in U\}$ is a face of $\SR(G)$, and the projection of $F$ onto the coordinates indexed by the nodes and edges of $G-U$ is an affine bijection onto $\SR(G-U)$. In particular, $\dist_{\SR(G)}(p,q)\le\diam(\SR(G-U))$ for any two vertices $p,q$ of $\SR(G)$ with $p_v=q_v=c_v$ for every $v\in U$.
\end{lemma}
\begin{proof}
The set $F$ is a face, since it is obtained by setting valid inequalities $z_v\ge0$ or $z_v\le1$ to equality. Let $z\in F$ and let $e=\{u,v\}\in E$ with $v\in U$. If $c_v=0$, then $0\le z_e\le z_v=0$, so $z_e=0$. If $c_v=1$, then $z_e\le z_u$ and $z_e\ge z_u+z_v-1=z_u$, so $z_e=z_u$. Conversely, every point of $\SR(G-U)$, extended by $z_v=c_v$ for $v\in U$ and by these formulas for the edges meeting $U$, satisfies all constraints of $\SR(G)$: if $z_e=0$ and $z_v=0$, the constraints of $e$ read $0\le z_u$ and $0\ge z_u-1$, and if $z_e=z_u$ and $z_v=1$, they read $0\le z_u\le 1$. Hence the projection is an affine bijection with affine inverse. Such a map preserves vertices and edges, and the one-skeleton of a face of $\SR(G)$ is a subgraph of the one-skeleton of $\SR(G)$.
\end{proof}

For $W\subseteq V$, we write $G[W]$ for the subgraph of $G$ induced by $W$. Let $p$ be a vertex of $\SR(G)$ and let $D_p:=\{v\in V:p_v\notin\{0,1\}\}$. By \cite[Theorem~7 and p.~157]{Pad89}, $p_v=\frac12$ for every $v\in D_p$. For an edge $e=\{u,v\}$ with $u,v\in D_p$, the constraints of $\SR(G)$ give $0\le p_e\le\frac12$. If $0<p_e<\frac12$, then $p$ is the midpoint of two points of $\SR(G)$ that differ only in the coordinate $e$. Hence $p_e\in\{0,\frac12\}$. We call $e$ \emph{odd} for $p$ if $p_e=0$ and \emph{even} for $p$ if $p_e=\frac12$. The rows involving $z_e$ that are tight at $p$ are $A^e$ and $C^e$ if $e$ is odd, and $B^e_u$ and $B^e_v$ if $e$ is even.

\begin{lemma}\label{lem:graph-round}
Let $p$ be a vertex of $\SR(G)$, let $K$ be the node set of a connected component of $G[D_p]$, and let $Y$ be a spanning tree of $G[K]$. Then $|K|\ge3$. Let $S\subseteq K$ satisfy $|S\cap e|=1$ for every edge $e$ of $Y$ that is odd for $p$, and $|S\cap e|\neq1$ for every edge $e$ of $Y$ that is even for $p$. There are exactly two such sets, and they are complementary in $K$. Let $p'$ be obtained from $p$ by setting $p'_v=1$ for $v\in S$, $p'_v=0$ for $v\in K\setminus S$, and $p'_e=\prod_{v\in e}p'_v$ for every edge $e$ meeting $K$. Then $p'$ is a vertex of $\SR(G)$ adjacent to $p$.
\end{lemma}
\begin{proof}
Along the edges of $Y$, the conditions on $S$ determine, for each node of $K$, whether it lies in $S$, once this is fixed for one node. This gives exactly two sets, complementary in $K$. The points $p$ and $p'$ agree on the nodes outside $D_p$, so by Lemma~\ref{lem:graph-face} we may work in $\SR(G[D_p])$. Since $G[D_p]$ is the disjoint union of its connected components, $\SR(G[D_p])$ is the Cartesian product of the polytopes $\SR(G[K'])$ over these components $K'$. The vertices of a product are the tuples of vertices of the factors, and two vertices that differ in exactly one factor are adjacent if and only if they are adjacent in that factor. Hence we may assume that $G$ is connected, $V=K=D_p$, and $p'=z^S$.

Since no node bound is tight at $p$ and the remaining edge constraints hold strictly, the constraints tight at $p$ are exactly the two listed above for each edge. As $p$ is a vertex, these $2|E|$ constraints have rank $|V|+|E|$, so $|E|\ge|V|$. Hence the connected graph $G$ contains a cycle, and $|V|\ge3$.

Call an edge $e$ \emph{satisfied} if $e$ is odd and $|S\cap e|=1$, or $e$ is even and $|S\cap e|\neq1$. All edges of $Y$ are satisfied. For a satisfied edge, both constraints tight at $p$ are tight at $z^S$. For an edge that is not satisfied, exactly one of them is: $C^e$ if $e$ is odd and $e\subseteq S$, $A^e$ if $e$ is odd and $e\cap S=\emptyset$, and $B^e_v$ if $e$ is even and $e\cap S=\{u\}$. In each case, one constraint tight at both points expresses $z_e$ as an affine function of the node coordinates. In addition, each satisfied edge contributes the equation $z_u+z_v=1$ if it is odd, and $z_u=z_v$ if it is even. Let $L$ be the affine space defined by the constraints tight at both $p$ and $z^S$. On $L$, every edge coordinate is an affine function of the node coordinates, and the node coordinates satisfy exactly these node equations. Hence the projection onto the node coordinates maps $L$ bijectively onto the solution space of these node equations, and the two have the same dimension. The associated homogeneous equations $w_u=-w_v$ and $w_u=w_v$ along the edges of the spanning tree $Y$ determine $w$ from its value at a single node, so $\dim L\le1$. The minimal face of $\SR(G)$ containing $p$ and $z^S$ lies in $L$ and contains two distinct points, so it is an edge of $\SR(G)$.
\end{proof}

\begin{proof}[Proof of Theorem~\ref{thm:diameter-sr-graph}]
If $E=\emptyset$, then $\SR(G)=[0,1]^V$, whose diameter is $|V|$. We prove the upper bound by induction on $|V|$. If $|V|=0$, then $\SR(G)$ is a point. Let $p\neq q$ be vertices of $\SR(G)$.

\textit{Case 1: some node $v$ has $p_v=q_v\in\{0,1\}$.} By Lemma~\ref{lem:graph-face} and induction, $\dist_{\SR(G)}(p,q)\le\diam(\SR(G-\{v\}))\le|V|-1$.

\textit{Case 2: $p$ and $q$ are integral.} Then $\dist_{\SR(G)}(p,q)\le|V|$ by Theorem~\ref{thm:diameter-mp}.

\textit{Case 3: there is $v\in D_p$ with $q_v\in\{0,1\}$.} Let $K$ be the node set of the connected component of $G[D_p]$ containing $v$, and choose the set $S$ of Lemma~\ref{lem:graph-round} with $v\in S$ if and only if $q_v=1$. The resulting vertex $p'$ is adjacent to $p$ and satisfies $p'_v=q_v$, so by Case~1 we have $\dist_{\SR(G)}(p,q)\le1+(|V|-1)=|V|$. The case where some $v\in D_q$ has $p_v\in\{0,1\}$ is symmetric.

\textit{Case 4: none of the above.} Then $D_p=D_q=:D$, and $D\neq\emptyset$. Let $K$ be the node set of a connected component of $G[D]$, let $Y$ be a spanning tree of $G[K]$, and let $S$ and $S'$ be sets as in Lemma~\ref{lem:graph-round} for $p$ and for $q$, respectively, with the same tree $Y$. For every edge $\{u,w\}$ of $Y$, the parities of $|S\cap\{u,w\}|$ and $|S'\cap\{u,w\}|$ coincide if and only if the edge is odd for both $p$ and $q$ or even for both. Hence, walking along $Y$, the indicator of whether a node belongs to exactly one of $S$ and $S'$ changes exactly at the edges of $Y$ whose type differs for $p$ and $q$. This partitions $K$ into two classes such that $S$ and $S'$ agree on one class and disagree on the other. Replacing $S'$ by $K\setminus S'$ exchanges the two roles. Since $|K|\ge3$, we can choose $S'$ so that $S$ and $S'$ agree on a class $W$ with $|W|\ge\lceil|K|/2\rceil\ge2$. Let $p'$ and $q'$ be the vertices given by Lemma~\ref{lem:graph-round}. They are adjacent to $p$ and $q$, respectively, and $p'_v=q'_v\in\{0,1\}$ for every $v\in W$. By Lemma~\ref{lem:graph-face} and induction, $\dist_{\SR(G)}(p',q')\le\diam(\SR(G-W))\le|V|-2$, so $\dist_{\SR(G)}(p,q)\le|V|$.
\end{proof}

\subsection{Large denominators}
\begin{proof}[Proof of Proposition~\ref{prop:denominators}]
Let $n=2\ell+1$ and $V=\{1,\dots,n\}$, and read indices modulo $n$, so that $n+1=1$. Let $E$ consist of $\{1,2\}$ and, for $r\in[\ell]$, of $p_r:=\{2r+1,2r+2\}$ and $t_r:=\{2r,2r+1,2r+2\}$. Then $H$ has rank three and $|E|=2\ell+1$. Consider the system
\[
z_1+z_2=1,\qquad z_{2r+1}=z_{2r+2},\qquad z_{2r}+z_{2r+1}+z_{2r+2}=2\qquad(r\in[\ell])
\]
in the node variables. Write $u_r:=z_{2r}$ for $r\in[\ell]$, and $u_{\ell+1}:=z_{2\ell+2}=z_1$. The last two equations give $u_{r+1}=1-\frac{u_r}2$ for $r\in[\ell]$, and hence $u_r-\frac23=(-\frac12)^{r-1}(u_1-\frac23)$ for $r\in[\ell+1]$. The first equation reads $u_1+u_{\ell+1}=1$. Therefore the system has a unique solution, given by
\[
u_1=\frac{2^\ell+2(-1)^\ell}{3\,(2^\ell+(-1)^\ell)},\qquad z_1=1-u_1=\frac{a}{2^\ell+(-1)^\ell},
\]
where $a:=(2^{\ell+1}+(-1)^\ell)/3$.
The number $a$ is an integer, since $2^{\ell+1}\equiv(-1)^{\ell+1}\pmod 3$. Any common divisor of $a$ and $2^\ell+(-1)^\ell$ divides $3a-2(2^\ell+(-1)^\ell)=(-1)^{\ell+1}$, so $z_1$ is in lowest terms. Since $\ell\ge2$, we have $0<u_1<1$. The map $u\mapsto1-\frac u2$ sends $(0,1)$ into $(\frac12,1)$, so all node values lie in $(0,1)$.

Complete these node values to a point $z\in\mathbb R^{V\cup E}$ by setting $z_{\{1,2\}}=0$, $z_{p_r}=z_{2r+1}$ and $z_{t_r}=0$ for $r\in[\ell]$. All inequalities of $\SR(H)$ hold. For $\{1,2\}$, we have $z_1+z_2-1=0$. For $p_r$, we have $z_{p_r}=z_{2r+1}=z_{2r+2}\ge 2z_{2r+1}-1$. For $t_r$, we have $z_{2r}+z_{2r+1}+z_{2r+2}-2=0$. Moreover, the following $2(2\ell+1)=|V|+|E|$ rows are tight at $z$: $A^e$ and $C^e$ for $e=\{1,2\}$ and for $e=t_r$, and $B^{p_r}_{2r+1}$ and $B^{p_r}_{2r+2}$, for $r\in[\ell]$. Use the rows $A^{\{1,2\}}$, $B^{p_r}_{2r+1}$ and $A^{t_r}$ to eliminate the edge variables. The remaining rows become exactly the system above, which has a unique solution. Hence these rows are linearly independent, and $z$ is a vertex of $\SR(H)$.
\end{proof}

\subsection{Low-edge elimination}
Let $\mathcal B$ be a basis at a vertex of $\SR(H)$. Recall that an edge $e$ is low if the rows of $\mathcal B$ involving $z_e$ are exactly $A^e$ and $C^e$. We call $e$ \emph{single} if exactly one row of $\mathcal B$ involves $z_e$.

\begin{lemma}\label{lem:fractional-node}
    Let $H$ be a hypergraph, and let $z\in\SR(H)$ satisfy $z_v\in\{0,1\}$ for every $v\in V$. Then $z$ is integral.
\end{lemma}

\begin{proof}
Let $e\in E$. If $z_v=0$ for some $v\in e$, then $0\le z_e\le z_v=0$, so $z_e=0$. Otherwise, $z_v=1$ for every $v\in e$, and $z_e\le 1$ together with $z_e\ge\sum_{v\in e}z_v-|e|+1=1$ gives $z_e=1$.
\end{proof}

\begin{lemma}\label{lem:active-system-integrality}
    Let $H$ be a hypergraph, let $z$ be a vertex of $\SR(H)$, and let $\mathcal{B}$ be a basis at $z$ with no low edge. Then $z$ is integral.
\end{lemma}

\begin{proof}
For $e\in E$, let $\mathcal A_e$ be the set of rows of $\mathcal B$ involving $z_e$.
Let $V_0$ be the set of nodes $v$ such that one of the following holds:
\begin{itemize}
\item[(i)] a node bound of $v$ is in $\mathcal B$,
\item[(ii)] $A^e,B^e_v\in\mathcal B$ for some $e\in E$,
\item[(iii)] $B^e_u,C^e\in\mathcal B$ for some $e\in E$ with $v\in e$ and some $u\in e\setminus\{v\}$.
\end{itemize}
We claim that $z_v\in\{0,1\}$ for every $v\in V_0$. In case (i) this is immediate.
In case (ii), $z_v=z_e=0$. In case (iii), $z_u=z_e=\sum_{r\in e}z_r-|e|+1$ gives
$\sum_{r\in e\setminus\{u\}}z_r=|e|-1$. Since $z_r\le 1$ for every $r$, all these
coordinates equal one, and in particular $z_v=1$.

Let $\Gamma$ be the graph on $V$ in which $u$ and $v$ are adjacent whenever
$B^e_u,B^e_v\in\mathcal B$ for some $e\in E$. For such $u,v$ we have $z_u=z_e=z_v$,
so $z$ is constant on each connected component of $\Gamma$. It therefore suffices
to show that every connected component $K$ of $\Gamma$ meets $V_0$: then $z_v\in\{0,1\}$
for every $v\in V$, and $z$ is integral by Lemma~\ref{lem:fractional-node}.

Suppose, for a contradiction, that some component $K$ is disjoint from $V_0$.
Define $w\in\mathbb R^{V\cup E}$ by $w_v=1$ for $v\in K$ and $w_v=0$ for $v\in V\setminus K$, and,
for every $e\in E$,
\[
w_e=\begin{cases}
w_u & \text{if } B^e_u\in\mathcal A_e \text{ for some } u\in e,\\
\sum_{v\in e}w_v & \text{if } \mathcal A_e=\{C^e\},\\
0 & \text{otherwise.}
\end{cases}
\]
The first case is well defined: if $B^e_u,B^e_{u'}\in\mathcal A_e$, then $u$ and $u'$
are adjacent in $\Gamma$, so $w_u=w_{u'}$. We show that $a^\top w=0$ for every row
$a^\top z=\beta$ of $\mathcal B$. Since $w\neq 0$, this contradicts the linear
independence of the rows of $\mathcal B$.
A node bound of $v$ in $\mathcal B$ gives $v\in V_0$, hence $v\notin K$ and $w_v=0$.
A row $B^e_u$ gives $w_e=w_u$ by definition. Consider a row $A^e$. If
$B^e_u\in\mathcal A_e$ for some $u$, then $u\in V_0$, so $w_e=w_u=0$. Otherwise $w_e=0$
by definition. Consider a row $C^e$. If $B^e_u\in\mathcal A_e$ for some $u$, then
$e\setminus\{u\}\subseteq V_0$, so $\sum_{v\in e}w_v=w_u=w_e$. If $\mathcal A_e=\{C^e\}$, the
equation holds by definition. In the remaining case $\mathcal A_e$ contains $C^e$,
contains no row $B^e_u$, and has at least two rows. Hence $\mathcal A_e=\{A^e,C^e\}$
and $e$ is a low edge, which is excluded by hypothesis.
\end{proof}

\begin{lemma}\label{lem:single-edge-preprocessing}
Let $\mathcal B$ be a basis at a vertex $z$ of $\SR(H)$, and let $s$ be the number of single edges of $\mathcal B$ whose row is $A^e$ or $C^e$. Then at most $s$ pivots transform $\mathcal B$ into a basis $\mathcal B'$ in which every single edge of $\mathcal B$ has a row $B^e_v$ with $v\in e$ as its only row, while all other rows of $\mathcal B$ are unchanged.
\end{lemma}

\begin{proof}
Let $e$ be a single edge of $\mathcal B$. If its row is $B^e_v$ for some $v\in e$, there is nothing to do. Otherwise, its row $r$ is $A^e$ or $C^e$. The other $d-1$ rows of $\mathcal B$ do not involve $z_e$ and have rank $d-1$, so they fix all coordinates except $z_e$. Along the line on which they are tight, the only constraints of $\SR(H)$ that vary are those of $e$, namely
\[
    \max\Bigl\{0,\sum_{v\in e}z_v-|e|+1\Bigr\}
    \le z_e\le
    \min_{v\in e}z_v.
\]
The row $r$ fixes $z_e$ at the lower end of this interval. Let $z'$ be the point of the line at which $z_e$ equals the upper end, and let $v\in e$ attain the minimum. Then $B^e_v$ is tight at $z'$, and replacing $r$ by $B^e_v$ gives $d$ linearly independent rows tight at $z'$. This is a pivot, which is degenerate if the two ends of the interval coincide. It changes only the row of $e$, and $e$ is still single. Repeating this for every single edge whose row is $A^e$ or $C^e$ uses at most $s$ pivots.
\end{proof}

\begin{proof}[Proof of Lemma~\ref{lem:frac-int}]
For every basis and every $e\in E$, let $\mathcal A_e$ be the set of rows of the basis that involve $z_e$. Every row other than a node bound involves exactly one edge coordinate. Since the rows of a basis have rank $d$, every edge coordinate is involved in some row, so $|\mathcal A_e|\ge1$. As the basis has $d=|V|+|E|$ rows, we obtain $\sum_{e\in E}(|\mathcal A_e|-1)\le|V|$. A low edge has $|\mathcal A_e|=2$, so every basis has at most $\min\{|V|,|E|\}$ low edges.

Let $s$ be the number of single edges of $\mathcal B$ whose row is $A^e$ or $C^e$, and let $b$ be the number of low edges of $\mathcal B$. Since a low edge is not single, we have $s+b\le|E|$. We first apply Lemma~\ref{lem:single-edge-preprocessing}, which uses at most $s$ pivots and does not change the rows of the other edges. Afterwards, every edge that is not low has some row $B^e_v$. Indeed, this holds for the single edges by construction, and an edge with at least two rows that has no row $B^e_v$ has exactly the rows $A^e$ and $C^e$, so it is low. We maintain this property while we eliminate the low edges one at a time.

Let $e$ be a low edge of the current basis $\mathcal B^\circ$ at the current vertex $z^\circ$. The rows of $\mathcal B^\circ$ other than $C^e$ have rank $d-1$ and define a line through $z^\circ$. Let $w$ be a direction of this line along which the slack of the inequality $z_e\ge\sum_{v\in e}z_v-|e|+1$ increases. Since $\SR(H)$ is bounded, there is a maximal $\lambda\ge0$ with $z^\circ+\lambda w\in\SR(H)$. Some inequality $q$ of $\SR(H)$ is tight at $z':=z^\circ+\lambda w$ and is violated by $z^\circ+\mu w$ for every $\mu>\lambda$. The other rows of $\mathcal B^\circ$ are tight along the whole line, so $q$ is none of them. Since the slack of $C^e$ increases along $w$, $q$ also differs from $C^e$. Since $q$ is not constant along $w$, its row is linearly independent of the rows of $\mathcal B^\circ$ other than $C^e$. Hence replacing $C^e$ by $q$ is a pivot to a basis at the vertex $z'$, which is degenerate if $\lambda=0$.

We now check the rows of each edge after this pivot. For every edge $f\neq e$, the rows of $f$ can only grow. If $f$ was not low, it keeps a row $B^f_v$. If $f$ was low and $q$ involves $z_f$, then $q$ is neither $A^f$ nor $C^f$, since both are already in the basis, so $q$ is a row $B^f_v$ and $f$ is no longer low. As for $e$, if $q$ involves $z_e$, then $q$ is a row $B^e_v$, since $A^e$ is still in the basis and $q\neq C^e$. Otherwise, $e$ is single with the row $A^e$. Every other single edge has a row $B^f_v$ by the property above, so Lemma~\ref{lem:single-edge-preprocessing} gives $e$ a row $B^e_v$ with one more pivot. In both cases, at most two pivots decrease the number of low edges by at least one and preserve the property above.

Since no low edge is created, after at most $b$ such eliminations we reach a basis $\hat{\mathcal B}$ with no low edge at a vertex $\hat z$. By Lemma~\ref{lem:active-system-integrality}, the vertex $\hat z$ is integral, so it is a multilinear vertex. The total number of pivots is at most
\[
s+2b\le (s+b)+b\le|E|+\min\{|V|,|E|\},
\]
where we used $b\le\min\{|V|,|E|\}$. Since every nondegenerate pivot moves along an edge of $\SR(H)$ and every degenerate pivot leaves the vertex unchanged, we also have $\dist_{\SR(H)}(z,\hat z)\le|E|+\min\{|V|,|E|\}$.
\end{proof}

\begin{proof}[Proof of Theorem~\ref{thm:diameter-sr}]
Let $p$ and $q$ be vertices of $\SR(H)$. Applying Lemma~\ref{lem:frac-int} to arbitrary bases at $p$ and $q$, we obtain multilinear vertices $\hat p$ and $\hat q$ with $\dist_{\SR(H)}(p,\hat p)\le\gamma$ and $\dist_{\SR(H)}(q,\hat q)\le\gamma$, where $\gamma:=|E|+\min\{|V|,|E|\}$. By Theorem~\ref{thm:diameter-mp}, we have $\dist_{\SR(H)}(\hat p,\hat q)\le|V|$, also if $\hat p=\hat q$. Hence
\begin{align*}
\dist_{\SR(H)}(p,q)&\le\dist_{\SR(H)}(p,\hat p)+\dist_{\SR(H)}(\hat p,\hat q)+\dist_{\SR(H)}(\hat q,q)\\
&\le|V|+2|E|+2\min\{|V|,|E|\}.
\end{align*}
\end{proof}

\section{Pivot Bounds}\label{app:pivots}
We use the basis exchange property of linear matroids, applied to the rows tight at a fixed vertex.
\begin{lemma}\label{lem:basis-exchange}
Let $\mathcal B$ and $\mathcal B'$ be bases at the same vertex $z$ of $\SR(H)$. Then $\mathcal B$ can be transformed into $\mathcal B'$ by $|\mathcal B'\setminus\mathcal B|\le d$ degenerate pivots. In particular, $\piv(\mathcal B,\mathcal B')\le d$.
\end{lemma}

\begin{proof}
Suppose that $\mathcal B\neq\mathcal B'$ and pick a row $r\in\mathcal B'\setminus\mathcal B$. Since $\mathcal B$ is a basis, the rows of $\mathcal B\cup\{r\}$ contain a unique minimal linearly dependent subset $\mathcal C$, and $r\in\mathcal C$. Since $\mathcal B'$ is linearly independent, $\mathcal C$ contains a row $r'\in\mathcal B\setminus\mathcal B'$. Then $(\mathcal B\setminus\{r'\})\cup\{r\}$ consists of $d$ linearly independent rows tight at $z$, so it is a basis at $z$. It is obtained from $\mathcal B$ by one pivot, which is degenerate since both bases are at $z$, and its intersection with $\mathcal B'$ is larger by one. Repeating this step $|\mathcal B'\setminus\mathcal B|$ times transforms $\mathcal B$ into $\mathcal B'$, and $|\mathcal B'\setminus\mathcal B|\le|\mathcal B'|=d$.
\end{proof}

\begin{proof}[Proof of Theorem~\ref{thm:pivot-flip}]
Let $S\setminus T=\{u_1,\dots,u_\mu\}$ and $T\setminus S=\{v_1,\dots,v_\nu\}$. Set $R_0:=S$, $R_i:=R_{i-1}\setminus\{u_i\}$ for $i\in[\mu]$, and $R_{\mu+j}:=R_{\mu+j-1}\cup\{v_j\}$ for $j\in[\nu]$. Then $R_{\mu+\nu}=T$, every $R_t$ contains $S\cap T$, and $R_t\subseteq S$ for $t\le \mu$, while $R_t\subseteq T$ for $t\ge \mu$. Moreover, $u_i\in R_t$ if and only if $t<i$, and $v_j\in R_t$ if and only if $t\ge \mu+j$.

Recall that $z^R_e=1$ if $e\subseteq R$ and $z^R_e=0$ otherwise. Hence $A^e$ is tight at $z^R$ if and only if $e\not\subseteq R$, and $B^e_w$ is tight at $z^R$ if and only if $e\subseteq R$ or $w\notin R$. For each $e\in E$ we define a row $\rho(e)$ as follows.
\begin{itemize}
\item If $e\subseteq S\cap T$, let $\rho(e):=B^e_w$ for an arbitrary $w\in e$. Since $e\subseteq S\cap T\subseteq R_t$, the row is tight at every $z^{R_t}$.
\item If $e\subseteq S$ and $e\not\subseteq T$, let $\rho(e):=B^e_{u_i}$, where $i$ is the smallest index with $u_i\in e$. For $t<i$, every node of $e\cap(S\setminus T)$ has index at least $i$ and still belongs to $R_t$, so $e\subseteq R_t$. For $t\ge i$, we have $u_i\notin R_t$. In both cases $\rho(e)$ is tight at $z^{R_t}$.
\item If $e\subseteq T$ and $e\not\subseteq S$, let $\rho(e):=B^e_{v_j}$, where $j$ is the largest index with $v_j\in e$. For $t\ge \mu+j$, every node of $e\cap(T\setminus S)$ has index at most $j$ and belongs to $R_t$, so $e\subseteq R_t$. For $t<\mu+j$, we have $v_j\notin R_t$. In both cases $\rho(e)$ is tight at $z^{R_t}$.
\item Otherwise, $e\not\subseteq S$ and $e\not\subseteq T$, and we let $\rho(e):=A^e$. Since every $R_t$ is contained in $S$ or in $T$, we have $e\not\subseteq R_t$, so $\rho(e)$ is tight at every $z^{R_t}$.
\end{itemize}
For $t\in\{0,\dots,\mu+\nu\}$, let $\mathcal B_t$ consist of the rows $\rho(e)$ for $e\in E$, together with $N^1_v$ for $v\in R_t$ and $N^0_v$ for $v\in V\setminus R_t$. These $d$ rows are tight at $z^{R_t}$. Setting these rows to equality fixes every node coordinate. Since $z_e$ appears in $\rho(e)$ with a nonzero coefficient, each $\rho(e)$ then fixes $z_e$. Hence these rows are linearly independent, and $\mathcal B_t$ is a basis at $z^{R_t}$. The bases $\mathcal B_{t-1}$ and $\mathcal B_t$ differ only in the node bound of the unique node of $R_{t-1}\triangle R_t$, so passing from one to the other is one nondegenerate pivot. With $\mathcal B^S:=\mathcal B_0$ and $\mathcal B^T:=\mathcal B_{\mu+\nu}$, the sequence $\mathcal B_0,\dots,\mathcal B_{\mu+\nu}$ consists of $\mu+\nu=|S\triangle T|$ pivots and visits only bases at the multilinear vertices $z^{R_t}$. Hence $\piv(\mathcal B^S,\mathcal B^T)\le|S\triangle T|\le|V|$.

If $E=\emptyset$, then $\SR(H)=[0,1]^V$. A degenerate pivot does not change the vertex, and a nondegenerate pivot moves to an adjacent vertex of the cube, which differs in one coordinate. Hence any bases at $z^S$ and $z^T$ are at least $|S\triangle T|$ pivots apart.

Finally, let $\mathcal B$ be any basis at $z^S$. By Lemma~\ref{lem:basis-exchange}, we have $\piv(\mathcal B,\mathcal B^S)\le d$, and
\[
\piv(\mathcal B,\mathcal B^T)\le\piv(\mathcal B,\mathcal B^S)+\piv(\mathcal B^S,\mathcal B^T)\le d+|S\triangle T|\le 2|V|+|E|.
\]
\end{proof}

\subsection{The component-switching construction}
Switching whole components instead of single nodes gives the following bound.
\begin{theorem}\label{thm:pivot-integral}
Let $H$ be a hypergraph and let $S,T\subseteq V$ be distinct. Let $C_1,\ldots,C_k$ be the connected components of $\Theta(H,S,T)$ in any order, let $S_0:=S$ and $S_i:=S_{i-1}\triangle C_i$ for $i\in[k]$, and let
\[
\mathcal M:=\{e\in E:e\not\subseteq S,\ e\not\subseteq T,\ e\subseteq S_i
\text{ for some }i\in[k-1]\}.
\]
Then there are bases $\mathcal B^S$ at $z^S$ and $\mathcal B^T$ at $z^T$ with $\piv(\mathcal B^S,\mathcal B^T)\le k+|\mathcal M|$. If $S,T$ satisfy $(\star)$, then $\mathcal M=\emptyset$, and hence $\piv(\mathcal B^S,\mathcal B^T)\le\co(\Theta(H,S,T))$.
\end{theorem}
\noindent Theorems~\ref{thm:pivot-flip} and~\ref{thm:pivot-integral} are incomparable. For some pairs $S,T$, the bound $k+|\mathcal M|$ exceeds $|S\triangle T|$ for every order of the components. For example, let $S\setminus T=\{a_1,a_1',a_2,a_2'\}$, $T\setminus S=\{b_1,b_2\}$ and $S\cap T=\{y_1,\dots,y_5\}$, and let $E$ consist of the edges $\{a_i,b_i\}$ and $\{a_i',b_i\}$ for $i\in[2]$, together with $\{a_1,a_1',b_2,y_j\}$ and $\{a_2,a_2',b_1,y_j\}$ for $j\in[5]$. Then $\Theta(H,S,T)$ has the two connected components $\{a_1,a_1',b_1\}$ and $\{a_2,a_2',b_2\}$, so $k=2$ and $|S\triangle T|=6$, while $|\mathcal M|=5$ for both orders. Under condition $(\star)$, the bound of Theorem~\ref{thm:pivot-integral} equals the distance bound $\co(\Theta(H,S,T))$ of Theorem~\ref{thm:diameter-mp}. This bound can be much smaller than $|S\triangle T|$.

We now describe the bases behind Theorem~\ref{thm:pivot-integral}. We construct them along the path of component switches, so that consecutive bases differ in few rows. This gives $k$ nondegenerate pivots and $|\mathcal M|$ degenerate pivots. Let $U:=S\triangle T$, let $C_1,\dots,C_k$ be the connected components of $\Theta(H,S,T)$, and let $S_0,\dots,S_k$ be as in Theorem~\ref{thm:pivot-integral}. We seek a basis at every $z^{S_i}$. Its rows of generating edges tie the node coordinates of each connected component to a representative node. Its node bounds then fix all node coordinates, and one selected row for each remaining edge fixes the corresponding edge coordinate.

For each component $C_l$, choose a representative node $r_l\in C_l$. Starting with $\Sigma:=\{r_l\}$, and while $\Sigma\neq C_l$, repeatedly choose a trace $g=f\cap U\in E_\Theta$, where $f$ is a generating edge, such that $g\cap\Sigma\neq\emptyset$ and $g\setminus\Sigma\neq\emptyset$. Choose $u\in g\cap\Sigma$. If $f$ satisfies condition (a) or (b) of Definition~\ref{def:contracted-hypergraph}, we select $B^f_u$ and $B^f_v$ for every $v\in g\setminus\Sigma$. If $f$ satisfies condition (c), then $g=\{u,v\}$, and we select $A^f$ and $C^f$. We then replace $\Sigma$ by $\Sigma\cup g$. Once an edge $f$ is selected, all nodes of $f\cap U$ are added to $\Sigma$, so $f$ is not selected again. Let $\mathcal F\subseteq F_\Theta$ be the set of generating edges selected in these constructions. We also select a node bound for each $v\in V\setminus U$ and for each representative node $r_l$, and, for each edge $e\in E\setminus\mathcal F$, one row among $A^e$ and $B^e_v$ with $v\in e$. Which of these rows is tight at $z^{S_i}$ depends on the set $I_e$ of Lemma~\ref{lem:pivot-interval}. Lemma~\ref{lem:adapted-basis} below shows that the selected rows form a basis at every $z^{S_i}$.

\begin{lemma}\label{lem:pivot-interval}
    For every $e\in E$, the set
    \[
        I_e:=\{i\in\{0,\ldots,k\}:e\subseteq S_i\}
    \]
    is either empty or an interval.
\end{lemma}

\begin{proof}
Each node of $U$ lies in exactly one component, so it changes sides exactly once along $S_0,\dots,S_k$, and the nodes outside $U$ never change sides. Hence, for each node $v$, the set $\{i:v\in S_i\}$ is empty or an interval. The set $I_e$ is the intersection of these intervals over $v\in e$, so it is empty or an interval as well.
\end{proof}

\begin{lemma}
\label{lem:adapted-basis}
Let $i\in\{0,\dots,k\}$. Consider the node bound tight at $z^{S_i}$ for every $v\in V\setminus U$ and for every representative node $r_l$, all rows selected above for the generating edges in $\mathcal F$, and, for every $e\in E\setminus\mathcal F$, one row among $A^e$ and $B^e_v$ with $v\in e$ that is tight at $z^{S_i}$. Then these $d$ rows form a basis at $z^{S_i}$.
\end{lemma}

\begin{proof}
We first verify that the selected generating edge rows are tight at $z^{S_i}$. If $f$ satisfies condition (a) or (b) of Definition~\ref{def:contracted-hypergraph}, then $f\cap U$ is contained in one component $C_l$. Hence all nodes of $f\cap U$ are switched simultaneously. Before that switch, $z_f$ and all $z_v$ with $v\in f\cap U$ have the same value, and the same holds after the switch. Thus every selected row $B^f_v$ is tight at every $z^{S_i}$. If $f$ satisfies condition (c), then
\[
    f\cap U = \{u,v\}, \quad u\in S\setminus T,\quad  v\in T\setminus S,
\]
and $u,v$ belong to the same component. Exactly one of $u,v$ therefore belongs to each $S_i$. Since $f\setminus U\subseteq S\cap T$, we have $z^{S_i}_f = 0$, and
\[
    \sum \limits_{x\in f} z_x^{S_i} = |f| - 1.
\]
Hence both $A^f$ and $C^f$ are tight at every $z^{S_i}$. We now count the selected rows. Each selected generating edge $f$ contributes one row for each node of $f\cap U$ newly added to $\Sigma$, plus one further row, namely $B^f_u$ or $A^f$. Over all $k$ component constructions, exactly $|U| - k$ non-representative nodes are added. Hence the generating edges contribute
\[
    |\mathcal{F}| + |U| - k
\]
rows. The selected node bounds contribute
\[
    |V| - |U| + k
\]
rows, and the edges in $E\setminus \mathcal{F}$ contribute
\[
    |E| - |\mathcal{F}|
\]
rows. In total there are 
\[
    |\mathcal{F}| + |U| - k + |V| - |U| + k + |E| - |\mathcal{F}| = |V| + |E|
\]
rows selected. It remains to prove that the selected rows are linearly independent. Consider the associated homogeneous system. The selected node bounds give 
\[
    w_v = 0, \quad v\in V\setminus U,
\]
and 
\[
    w_{r_l} = 0, \quad l\in [k].
\]
For every generating edge $f\in \mathcal{F}$ that satisfies condition (a) or (b), the selected rows imply 
\[
    w_u = w_v
\]
for the nodes of $f\cap U$ linked by construction. For a generating edge $f$ that satisfies condition (c), the selected rows $A^f$ and $C^f$, together with $w_x = 0$ for $x\in f\setminus U$, imply
\[
    w_u + w_v = 0, \quad f\cap U = \{u,v\}.
\]
By construction, these relations connect every node of $C_l$ to its representative node $r_l$. Hence every $v\in C_l$ satisfies
\[
    w_v = \pm w_{r_l} = 0.
\]
Since $C_1,\dots, C_k$ partition $U$, it follows that 
\[
    w_v = 0, \quad v\in V.
\]
Now let $f\in \mathcal{F}$. If $f$ satisfies condition (a) or (b), any selected row $B^f_v$ gives
\[
    w_f = w_v = 0.
\]
If $f$ satisfies condition (c), the row $A^f$ gives $w_f = 0$. Thus
\[
    w_f = 0,\quad f\in \mathcal{F}.
\]
Finally, for every $e\in E\setminus \mathcal{F}$, the selected row is either $A^e$, which gives $w_e = 0$, or some $B^e_v$, which gives
\[
    w_e = w_v = 0.
\]
Therefore $w = 0$. Hence the selected $d$ rows are linearly independent and form a basis at $z^{S_i}$.
\end{proof}

\begin{lemma}\label{lem:middle-star}
Let $S,T\subseteq V$ be distinct, and let $C_1,\dots,C_k$, $S_0,\dots,S_k$, and $\mathcal M$ be as in Theorem~\ref{thm:pivot-integral}. Then every edge $e\in\mathcal M$ violates condition $(\star)$. In particular, if $S,T$ satisfy $(\star)$, then $\mathcal M=\emptyset$.
\end{lemma}
\begin{proof}
Let $e\in\mathcal M$. Since $e\subseteq S_i\subseteq S\cup T$ for some $i$, $e\not\subseteq S$ and $e\not\subseteq T$, the edge $e$ is contained in $S\cup T$ and meets both $S\setminus T$ and $T\setminus S$. Suppose that $e$ contains $u\in S\setminus T$ and $v\in T\setminus S$ in the same component $C_l$. If $l\le i$, then $u\notin S_i$, and if $l>i$, then $v\notin S_i$. In both cases $e\not\subseteq S_i$, a contradiction.
\end{proof}

\begin{proof}[Proof of Theorem~\ref{thm:pivot-integral}]
We show that every edge $e\in E\setminus\mathcal F$ has a row that is tight at every $z^{S_i}$, except when $I_e=\{\alpha,\dots,\beta\}$ with $0<\alpha\le\beta<k$. For such an edge one degenerate pivot exchanges its row, and these edges are exactly the edges of $\mathcal M$. Recall that $\mathcal{F}$ is the set of generating edges selected in the construction preceding Lemma~\ref{lem:adapted-basis}. For each $e\in E\setminus \mathcal{F}$, $I_e$ is either empty or an interval of consecutive integers by Lemma~\ref{lem:pivot-interval}. If $I_e = \emptyset$, then $z^{S_i}_e = 0$ for every $i$, so $A^e$ is tight throughout the path. If $I_e = \{0,\dots, k\}$, then $e\subseteq S_i$ for every $i$. Hence, for any $u\in e$, we have
\[  
    z^{S_i}_e = z^{S_i}_u = 1
\]
for every $i$, so $B^e_u$ is tight throughout. Now suppose $I_e = \{0,\dots, \beta\}, \beta< k$. Since
\[
    e\subseteq S_\beta,\quad e\nsubseteq S_{\beta+1},
\]
there exists $v\in e\cap C_{\beta + 1}$ whose node coordinate changes from one to zero at step $\beta + 1$. Therefore
\[
    z^{S_i}_e = z^{S_i}_v = \begin{cases}
        1, \quad 0\le i\le \beta,\\
        0, \quad \beta + 1\le i\le k,
    \end{cases}
\]
and $B^e_v$ is tight throughout the path. Similarly, suppose 
\[
    I_e = \{\alpha, \dots, k\},\quad \alpha > 0.
\]
Since
\[
    e\nsubseteq S_{\alpha - 1}, e\subseteq S_\alpha,
\]
there exists $u\in e\cap C_\alpha$ whose node coordinate changes from zero to one at step $\alpha$. Hence
\[
    z^{S_i}_e = z^{S_i}_u = \begin{cases}
        0, \quad 0\le i < \alpha,\\
        1, \quad \alpha \le i\le k,
    \end{cases}
\]
so $B^e_u$ is tight throughout. It then remains to consider 
\[
    I_e = \{\alpha, \dots, \beta\}, \quad 0<\alpha \le \beta < k.
\]
Since $e\nsubseteq S_{\alpha - 1}$ but $e\subseteq S_\alpha$, there exists $u\in e\cap C_\alpha$ whose node coordinate changes from zero to one at step $\alpha$. Similarly, since $e\subseteq S_\beta$ but $e\nsubseteq S_{\beta + 1}$, there exists $v\in e\cap C_{\beta + 1}$ whose node coordinate changes from one to zero at step $\beta + 1$. Consequently, 
\[
    z^{S_i}_e =\begin{cases}
        z^{S_i}_u, \quad 0\le i\le \beta,\\
        z^{S_i}_v, \quad \alpha \le i \le k.
    \end{cases}
\]
We use $B^e_u$ before reaching $z^{S_\alpha}$, and at $z^{S_\alpha}$ replace it by $B^e_v$. By Lemma~\ref{lem:adapted-basis}, either choice, together with the other selected rows, forms a basis at $z^{S_\alpha}$. Hence this replacement is one degenerate pivot. Moreover, if $f\in \mathcal{F}$ satisfies condition (a), then $f\subseteq S$, and if it satisfies condition (b), then $f\subseteq T$. If it satisfies condition (c), the two nodes of $f\cap U$ lie on opposite sides of $S\triangle T$ and in the same component, so they are switched simultaneously. Hence $f\nsubseteq S_i$ for every $i$. Therefore, $\mathcal{F}\cap \mathcal{M} = \emptyset$, and the last case occurs precisely for the edges in $\mathcal{M}$. Thus the above construction requires exactly $|\mathcal{M}|$ degenerate pivots. We now traverse the path of component switches. Consider the transition
\[
    z^{S_{i-1}} \to z^{S_i}.
\]
Before making any degenerate exchange scheduled at $z^{S_i}$, the row selected for every $e\in E\setminus \mathcal{F}$ is tight at both endpoints. The rows selected from the generating edges in $\mathcal{F}$ are also tight at both endpoints. Among the selected node bounds, only the node bound of the representative $r_i\in C_i$ changes. Therefore, replacing the node bound of $r_i$ that is tight at $z^{S_{i-1}}$ by the node bound that is tight at $z^{S_i}$ changes exactly one row. By Lemma~\ref{lem:adapted-basis}, the resulting rows form a basis at $z^{S_i}$. Hence this replacement is one nondegenerate pivot. After reaching $z^{S_i}$, we perform the degenerate exchanges corresponding to the edges $e\in\mathcal{M}$ whose interval $I_e$ begins at $i$. Repeating this procedure for $i\in[k]$ yields bases $\mathcal{B}^S$ at $z^S$ and $\mathcal{B}^T$ at $z^T$ connected by $k$ nondegenerate pivots and $|\mathcal{M}|$ degenerate pivots. Therefore $\piv(\mathcal{B}^S,\mathcal{B}^T)\le k+|\mathcal{M}|$. Finally, if $S,T$ satisfy $(\star)$, then $\mathcal M=\emptyset$ by Lemma~\ref{lem:middle-star}, and $k=\co(\Theta(H,S,T))$.
\end{proof}

\begin{proof}[Proof of Proposition~\ref{prop:lowerbound}]
    Let $V = \{a,b_1,\dots, b_n\}$, and define
    \[
        E = \{e\subseteq V: a\in e, |e|\ge 2\}.
    \]
    Then $|V| = n+1$ and $|E| = 2^n - 1$. Define
    \[
        S := V, \quad T:= V\setminus \{a\}.
    \]
    Since $S\triangle T = \{a\}$, the contracted hypergraph $\Theta(H,S,T)$ is connected. By Theorem~\ref{thm:adjacency-sr}, $z^S$ and $z^T$ are adjacent in $\SR(H)$. We first construct a basis $\mathcal{B}^S$ at $z^S$. For every $e\in E$, choose an arbitrary node $b(e)\in e\setminus \{a\}$, and define
    \[
        \mathcal{B}^S := \{N^1_v: v\in V\} \cup \{B^e_{b(e)}: e\in E\}.
    \]
    All these rows are tight at $z^S$, since every node and edge coordinate of $z^S$ equals one. Moreover, the node bounds determine all node coordinates, while each selected row $B^e_{b(e)}$ determines the corresponding edge coordinate $z_e$. Thus the $|V| + |E|$ selected rows are linearly independent and form a basis at $z^S$. We now prove the lower bound. At $z^T$, we have 
    \[
        z^T_a = 0, \quad z^T_{b_i} = 1, \forall i\in [n].
    \]
    Since every edge $e\in E$ contains $a$, we also have 
    \[
        z^T_e = 0, \forall e\in E.
    \]
    Consequently, none of the rows $B^e_{b(e)}$ selected in $\mathcal{B}^S$ is tight at $z^T$, because
    \[
        z^T_e = 0 \neq 1 = z^T_{b(e)}.
    \]
    Furthermore, the row $N^1_a$ is not tight at $z^T$. Therefore, every basis $\mathcal{B}^T$ at $z^T$ excludes at least $|E| + 1$ rows of $\mathcal{B}^S$. Hence 
    \[
        |\mathcal{B}^S\setminus \mathcal{B}^T| \ge |E| + 1.
    \]
    Since each pivot replaces exactly one row, we conclude that 
    \[
        \piv(\mathcal{B}^S,\mathcal{B}^T) \ge |E| + 1
    \]
    for every basis $\mathcal{B}^T$ at $z^T$. It remains to show that equality is attainable. Starting from $\mathcal{B}^S$, for every $e\in E$, replace the row $B^e_{b(e)}$ by $B^e_a$. Both rows are tight at $z^S$. Since all node bounds remain selected and exactly one row involving each edge coordinate is retained, each replacement preserves linear independence. Hence each replacement is a degenerate pivot. After $|E|$ such pivots, we obtain the basis
    \[
        \hat{\mathcal{B}}^S = \{N^1_v: v\in V\}\cup \{B^e_a: e\in E\}.
    \]
    Now replace the row $N^1_a$ by $N^0_a$, obtaining 
    \[
        \mathcal{B}^T = \{N^0_a\}\cup \{N^1_{b_i}: i\in [n]\}\cup \{B^e_a: e\in E\}.
    \]
    All these rows are tight at $z^T$. They are also linearly independent, since the node bounds determine all node coordinates and each row $B^e_a$ determines the corresponding edge coordinate. Thus $\mathcal{B}^T$ is a basis at $z^T$. To verify that this replacement corresponds to a nondegenerate pivot, consider the face defined by the $|V| + |E| - 1$ common rows:
    \[
        F = \{z\in \SR(H): z_{b_i} = 1\, \forall i\in [n], z_e = z_a \, \forall e\in E\}.
    \]
    Writing $t = z_a$, the inequalities defining $\SR(H)$ imply $0\le t\le 1$. Conversely, every $t\in [0,1]$ defines a feasible point of $F$ by setting $z_{b_i} = 1$ for every $i\in [n]$ and $z_e = t$ for every $e\in E$. Hence $F$ is a one-dimensional face whose endpoints are precisely $z^S$ and $z^T$. Therefore, replacing $N^1_a$ by $N^0_a$ is one nondegenerate pivot from $z^S$ to $z^T$. We have now constructed a sequence of $|E|$ degenerate pivots followed by one nondegenerate pivot, giving
    \[
        \piv(\mathcal{B}^S, \mathcal{B}^T) \le |E| + 1.
    \]
    Combining this with the lower bound, we obtain
    \[
        \piv(\mathcal{B}^S, \mathcal{B}^T) = |E| + 1 = 2^n.
    \]
\end{proof}

\begin{proof}[Proof of Theorem~\ref{thm:pivot-sr}]
Let $\gamma:=|E|+\min\{|V|,|E|\}$. By Lemma~\ref{lem:frac-int}, there exist multilinear vertices $z^S$ and $z^T$ and bases $\hat{\mathcal B}$ at $z^S$ and $\hat{\mathcal B}'$ at $z^T$ such that $\piv(\mathcal B,\hat{\mathcal B})\le\gamma$ and $\piv(\mathcal B',\hat{\mathcal B}')\le\gamma$. Since pivots are reversible, we also have $\piv(\hat{\mathcal B}',\mathcal B')\le\gamma$. By Theorem~\ref{thm:pivot-flip}, there exist bases $\mathcal B^S$ at $z^S$ and $\mathcal B^T$ at $z^T$ such that $\piv(\mathcal B^S,\mathcal B^T)\le|S\triangle T|\le|V|$. By Lemma~\ref{lem:basis-exchange}, we have $\piv(\hat{\mathcal B},\mathcal B^S)\le d$ and $\piv(\mathcal B^T,\hat{\mathcal B}')\le d$. Therefore
\begin{align*}
    \piv(\mathcal{B},\mathcal{B}')
    &\le \piv(\mathcal{B},\hat{\mathcal{B}})
       + \piv(\hat{\mathcal{B}},\mathcal{B}^S)
       + \piv(\mathcal{B}^S,\mathcal{B}^T)
       + \piv(\mathcal{B}^T,\hat{\mathcal{B}}')
       + \piv(\hat{\mathcal{B}}',\mathcal{B}') \\
    &\le 2\gamma+2d+|V| = 3|V|+4|E|+2\min\{|V|,|E|\}.
\end{align*}
\end{proof}
\end{document}